\documentclass{article}

\usepackage[a4paper]{geometry}

\usepackage[francais,english]{babel}
\usepackage{amsfonts,latexsym,amstext}
\usepackage{amsmath,amscd,amssymb}

\usepackage[latin1]{inputenc}

\usepackage{aeguill}

\usepackage[T1]{fontenc}

\usepackage[usenames]{color}

\usepackage[a4paper]{geometry}

\usepackage{url}

\usepackage{psfrag}

\usepackage{placeins}

\usepackage[dvips]{graphicx}
\usepackage{latexsym,amssymb,amsmath,color,times,bbm}

\newtheorem{thm}{Theorem}[section]

\newtheorem{lem}[thm]{Lemma}
\newtheorem{prop}[thm]{Proposition}
\newtheorem{rem}[thm]{Remark}

\numberwithin{equation}{section}
\newcommand{\norm}[1]{\left\Vert#1\right\Vert}
\newcommand{\abs}[1]{\left\vert#1\right\vert}
\newcommand{\set}[1]{\left\{#1\right\}}

\def\R{\mathbb R}
\def\N{\mathbb N}

\def\E{\mathbb E}
\def\P{\mathbb P}
\def\Q{\mathbb Q}

\def\cprime{$'$}

\newcommand{\tr}[1]{{\vphantom{#1}}^{\mathit t}{#1}}

\def\sqw{\hbox{\rlap{\leavevmode\raise.3ex\hbox{$\sqcap$}}$%
\sqcup$}}
\def\sqb{\hbox{\hskip5pt\vrule width4pt height6pt depth1.5pt%
\hskip1pt}}

\def\qed{\ifmmode\hbox{\hfill\sqb}\else{\ifhmode\unskip\fi%
\nobreak\hfil
\penalty50\hskip1em\null\nobreak\hfil\sqb
\parfillskip=0pt\finalhyphendemerits=0\endgraf}\fi}
\def\cqfd{\ifmmode\sqw\else{\ifhmode\unskip\fi\nobreak\hfil
\penalty50\hskip1em\null\nobreak\hfil\sqw
\parfillskip=0pt\finalhyphendemerits=0\endgraf}\fi}

\begin{document}

\renewcommand{\labelitemi}{$\bullet$}
\bibliographystyle{plain}
\pagestyle{headings}
\title{Markovian quadratic and superquadratic BSDEs with an unbounded terminal condition}

\author{
Adrien Richou\\
Univ. Bordeaux, IMB, UMR 5251, F-33400 Talence, France.\\
CNRS, IMB, UMR 5251, F-33400 Talence, France.\\
INRIA, \'Equipe ALEA, F-33400 Talence, France.\\ 
e-mail: adrien.richou@math.u-bordeaux1.fr}

\selectlanguage{english}

\maketitle

\begin{abstract} 
This article deals with the existence and the uniqueness of solutions to quadratic and superquadratic Markovian backward stochastic differential equations (BSDEs for short) with an unbounded terminal condition. Our results are deeply linked with a strong a priori estimate on $Z$ that takes advantage of the Markovian framework. This estimate allows us to prove the existence of a viscosity solution to a semilinear parabolic partial differential equation with nonlinearity having quadratic or superquadratic growth in the gradient of the solution. This estimate also allows us to give explicit convergence rates for time approximation of quadratic or superquadratic Markovian BSDEs. 
\end{abstract}

%
%
\selectlanguage{english}



\section{Introduction}

Since the early nineties and the work of Pardoux and Peng \cite{Pardoux-Peng-90}, there has been an increasing interest for backward stochastic differential equations (BSDEs for short). These equations have a wide range of applications in stochastic control, in finance or in partial differential equation theory. A particular class of BSDE is studied since few years: BSDEs with generators of quadratic growth with respect to the variable $z$ (quadratic BSDEs for short). This class arises, for example, in the context of utility optimization problems with exponential utility functions, or alternatively in questions related to risk minimization for the entropic risk measure (see e.g. \cite{Rouge-ElKaroui-00,Hu-Imkeller-Muller-05,Mania-Schweizer-05} among many other references). Many papers deal with existence and uniqueness of solution for such BSDEs. In the first one \cite{Kobylanski-00}, Kobylanski obtains an existence and uniqueness result for quadratic BSDEs when the terminal condition is bounded. Let us remark that this result has been revisited recently thanks to a fixed point argument by Tevzadze in \cite{Tevzadze-08}. Now, it is well known that the boundedness of the terminal condition is a too strong assumption. Indeed, when we look to the simple quadratic BSDE
$$Y_t=\xi+\int_t^T \frac{\abs{Z_s}^2}{2}ds-\int_t^T Z_s dW_s,$$
we find the explicit solution $Y_t=\log \left(\E\left[e^{\xi}|\mathcal{F}_t \right] \right)$ and we immediately see that we just need to have an exponential moment for $\xi$ to obtain a solution. In \cite{Briand-Hu-06}, Briand and Hu show an existence result for quadratic BSDEs when the terminal condition has such an assumption. Let us notice that this result has been recently revisited in \cite{Barrieu-ElKaroui-11} by a direct forward method that does not use the result of Kobylanski. In this paper, Barrieu and El Karoui obtain a monotone stability result for general quadratic semimartingales and then derive an existence result for general quadratic BSDEs. For the uniqueness problem, results are more incomplete. In \cite{Delbaen-Hu-Richou-09}, authors show a uniqueness result when the generator is convex (or concave) with respect to $z$ and when $\xi$ has an exponential moment which is almost the exponential moment needed for the existence result.

Naturally, we could also wonder what happens when the generator has a superquadratic growth with respect to the variable $z$. Up to our knowledge the case of superquadratic BSDEs is only investigate in the recent paper \cite{Delbaen-Hu-Bao-09}. In this article, authors consider superquadratic BSDEs when the terminal condition is bounded and the generator is convex in $z$. Firstly, they show that in a general way the problem is ill-posed: given a superquadratic generator, there exists a bounded terminal condition such that the associated BSDE does not admit any bounded solution and, on the other hand, if the BSDE admits a bounded solution, there exist infinitely many bounded solutions for this BSDE. In the same paper, authors also show that the problem becomes well-posed in a Markovian framework: When the terminal condition and the generator are deterministic functions of a forward SDE, we have an existence result.

The first aim of this paper is to study existence and uniqueness results for quadratic and superquadratic Markovian BSDEs. More precisely, we consider $(X,Y,Z)$ the solution to the (decoupled) forward backward system
\begin{eqnarray*}
X_t &=& x +\int_0^t b(s,X_s)ds+\int_0^t \sigma(s) dW_s,\\
Y_t &=& g(X_T) +\int_t^T f(s,X_s,Y_s,Z_s)ds-\int_t^T Z_s dW_s,
\end{eqnarray*}
where $f$ has quadratic or superquadratic growth with respect $z$, has no convexity assumption, and $g$ is not supposed to be bounded. The starting point of our work is a simple result that says: if $g$ and $f$ are Lipschitz functions with respect to $x$, then there exists a unique solution such that $Z$ is bounded, or in other words, $Z$ preserves the regularity of the derivatives of $g$ and $f$ with respect to $x$. Now, the idea is to show that this property stays true when $g$ and $f$ are only locally Lipschitz. More precisely, if we assume that 
$$\abs{\nabla g(x)} + \abs{\nabla_x f(.,x,.,.)} \leqslant C(1+\abs{x}^r)$$
for $r$ sufficiently small, then we are able to show the a priori estimate
$$\abs{Z} \leqslant C(1+\abs{X}^r).$$ 
Thanks to this kind of estimate, it is then possible to show an existence and uniqueness result amongst solutions that, roughly speaking, verify such an estimate (see Theorem~\ref{theoreme r=1surl}). Contrarily to \cite{Delbaen-Hu-Richou-09,Delbaen-Hu-Bao-09} we do not need a convexity assumption on $f$ and contrarily to \cite{Delbaen-Hu-Bao-09}, we treat the case of unbounded terminal conditions. On the other hand, for the quadratic case we need a framework which is far more restrictive than the general framework of \cite{Delbaen-Hu-Richou-09,Barrieu-ElKaroui-11} because we only consider markovian BSDEs with assumptions on the derivatives of $g$ and $f$ instead of classical assumptions on the growth of $g$ and $f$.

One of the major drawback of results explained before is that we consider only the case of a deterministic function $\sigma$. The second part of our paper gives some partial results when $\sigma$ is random. In this framework we do not know if our previous starting point stays true: if $g$ and $f$ are Lipschitz functions with respect to $x$ and if $\sigma$ is bounded, does there exist a solution such that $Z$ is bounded ? We are able to show that this is true when $T$ is small enough or for all $T$ when we consider a simple example of quadratic BSDE. But the general case stays an open question. We also investigate precisely the quadratic case when $g$ and $f$ are bounded with respect to $x$ by deeply using bounded mean oscillation martingale (BMO martingale for short) tools. 

Thanks to our existence and uniqueness result we are able to give a probabilistic representation of the following PDE:
\begin{equation*}
\left\{
\begin{array}{l}
\partial_t u(t,x)+\mathcal{L}u(t,x)+f(t,x,u(t,x),^t\nabla u(t,x)\sigma(t))=0, \quad x \in \R^d, t \in [0,T],\\
u(T,.)=g.
\end{array}
\right.
\end{equation*}
Such a probabilistic representation, also called Feynman-Kac representation, is already given in \cite{Delbaen-Hu-Richou-09} when $f$ has a quadratic growth and is convex with respect to $z$. Existence and uniqueness of this PDE has been studied in \cite{DaLio-Ley-06} when $f$ has a quadratic growth with respect to $\nabla u\sigma$ and in \cite{DaLio-Ley-08} for the superquadratic case, but the main part of the results needs a convexity assumption on $f$ with respect to $z$. In this paper, our existence result arises in quadratic and superquadratic frameworks. Moreover, we do not need any convexity assumption on $f$ which is interesting for applications: For example, when we consider Isaacs equations in differential game theory, $f$ is the sum of a convex function and a concave function with respect to $z$.

The main goal of this paper is to apply a priori estimates obtained for the process $Z$ to the problem of time discretization of quadratic and superquadratic BSDEs. Actually, the design of efficient algorithms which are able to solve BSDEs in any reasonable dimension has been intensively studied since the first work of Chevance~\cite{Chevance-97}, see for instance \cite{Zhang-04,Bouchard-Touzi-04,Gobet-Lemor-Warin-05}. But in all these works, the driver of the BSDE is a Lipschitz function with respect to $z$ and this assumption plays a key role in theirs proofs. In a recent paper, Cheridito and Stadje~\cite{Cheridito-Stadje-10} study approximation of BSDEs by backward stochastic difference equations which consist in replacing the Brownian motion by a random walk. They obtain a convergence result when the driver has a subquadratic growth with respect to $z$ and they give an example where there proof does not work when the driver has a quadratic growth. To the best of our knowledge, the only works where the time approximation of a quadratic BSDE is studied are the one of Imkeller and dos Reis~\cite{Imkeller-dosReis-09} and the one of Richou~\cite{Richou-11}. Let us notice that, when the driver has a specific form\footnote{Roughly speaking, the driver is a sum of a quadratic term $z \mapsto C\abs{z}^2$ and a function that has a linear growth with respect to $z$.}, it is possible to get around the problem by using an exponential transformation method (see~\cite{Imkeller-Reis-Zhang-10}) or by using results on fully coupled forward-backward differential equations (see~\cite{Delarue-Menozzi-06}). Papers \cite{Imkeller-dosReis-09,Richou-11} only study the case of a bounded terminal condition: The first one investigates the case of Lipschitz terminal conditions whereas the second one studies the non-smooth case. To the best of our knowledge, the time approximation of superquadratic BSDEs has not been studied yet. In this paper we have obtained two types of results. Firstly we consider the case of a deterministic function $\sigma$. Theorem~\ref{theorem vitesse convergence sigma deterministe} gives us a speed of convergence very close to the speed of convergence in the classical Lipschitz case and this theorem is obtained in a general framework (quadratic and superquadratic BSDEs with an unbounded terminal condition). When $\sigma$ is random, we only study quadratic BSDEs with bounded terminal conditions. In Theorem~\ref{theorem vitesse convergence sigma aleatoire} we obtain almost the classical speed of convergence but in a restricted framework that does not cover some interesting situations: for example we are not able to find a ``good'' speed of convergence when $\sigma$ and $g$ are Lipschitz functions with respect to $x$ and this question is actually a real challenge.

The paper is organized as follows. In section 2 we obtain an existence and uniqueness result and an a priori estimate on $Z$ for Markovian quadratic and superquadratic BSDEs with unbounded terminal conditions when $\sigma$ is a deterministic function. In section 3 we give some extra partial results when $\sigma$ is random. Section 4 contains an application to semilinear parabolic PDEs. The last section is devoted to time approximation of quadratic and superquadratic Markovian BSDEs.

\paragraph{Notations}
Throughout this paper, $(W_t)_{t \geqslant 0}$ will denote a $d$-dimensional Brownian motion, defined on a probability space $(\Omega,\mathcal{F}, \P)$. For $t \geqslant 0$, let $\mathcal{F}_t$ denotes the $\sigma$-algebra $\sigma(W_s; 0\leqslant s\leqslant t)$, augmented with the $\P$-null sets of $\mathcal{F}$. The Euclidean norm on $\R^d$ will be denoted by $|.|$. The operator norm induced by $|.|$ on the space of linear operators is also denoted by $|.|$. The notation $\mathbb{E}_t$ stands for the conditional expectation given $\mathcal{F}_t$. For $p \geqslant 2$, $m \in \N$, we denote further
\begin{itemize}
 \item $\mathcal{S}^p(\R^m)$, or $\mathcal{S}^p$ when no confusion is possible, the space of all adapted processes $(Y_t)_{t \in [0,T]}$ with values in $\R^m$ normed by $\norm{Y}_{\mathcal{S}^p}=\E [(\sup_{t \in [0,T]} \abs{Y_t})^p]^{1/p}$; $\mathcal{S}^{\infty}(\R^m)$, or $\mathcal{S}^{\infty}$, the space of bounded measurable processes;
 \item $\mathcal{M}^p(\R^m)$, or $\mathcal{M}^p$, the space of all progressively measurable processes $(Z_t)_{t \in [0,T]}$ with values in $\R^m$ normed by $\norm{Z}_{\mathcal{M}^p}=\E[(\int_0^T \abs{Z_s}^2ds)^{p/2}]^{1/p}$.
\end{itemize}
 In the following, we keep the same notation $C$ for all finite, nonnegative constants that appear in our computations.

In this paper we will consider $X$ the solution to the SDE
\begin{equation}
\label{EDS}
 X_t=x+\int_0^t b(s,X_s)ds+\int_0^t \sigma(s) dW_s,
\end{equation}
and $(Y,Z) \in \mathcal{S}^2\times \mathcal{M}^2$ the solution to the Markovian BSDE
\begin{equation}
\label{EDSR}
 Y_t=g(X_T)+\int_t^T f(s,X_s,Y_s,Z_s)ds-\int_t^T Z_sdW_s.
\end{equation}

\section{A uniqueness and existence result}
For the SDE (\ref{EDS}) we use standard assumption.
\paragraph{Assumption (F.1).}
Let $b : [0,T] \times \mathbb{R}^d \rightarrow \mathbb{R}^d$ and $\sigma : [0,T] \rightarrow \mathbb{R}^{d \times d}$ be continuous  functions and let us assume that there exists  $K_b \geqslant 0$ such that:
\begin{enumerate}
 \item $\forall t \in [0,T]$, $\abs{b(t,0)} \leqslant C$,
 \item $\forall t \in [0,T]$, $\forall (x,x') \in \mathbb{R}^d \times \mathbb{R}^d$, $\abs{b(t,x)-b(t,x')} \leqslant K_b \abs{x-x'}.$
\end{enumerate}
Let us assume the following for the generator and the terminal condition of the BSDE (\ref{EDSR}).
\paragraph{Assumption (B.1).}
Let $f: [0,T] \times \mathbb{R}^d \times \mathbb{R} \times \mathbb{R}^{1\times d} \rightarrow \mathbb{R}$ and $g:\mathbb{R}^d \rightarrow  \mathbb{R}$ be continuous functions and let us assume moreover that there exist five constants, $l \geqslant 1$, $\alpha \geqslant 0$, $\beta \geqslant 0$, $\gamma \geqslant 0$ and $K_{f,y} \geqslant 0$ such that:
\begin{enumerate}
\item for each $(t,x,y,y',z) \in [0,T] \times \mathbb{R}^d \times \mathbb{R} \times \mathbb{R} \times \mathbb{R}^{1\times d}$,
$$ \abs{f(t,x,y,z)-f(t,x,y',z)} \leqslant K_{f,y} \abs{y-y'};$$
\item for each $(t,x,y,z,z') \in [0,T] \times \R^d \times \R \times \R^{1\times d} \times \R^{1\times d}$,
$$\abs{ f(t,x,y,z)-f(t,x,y,z')} \leqslant \left(C+\frac{\gamma}{2}(\abs{z}^l+\abs{z'}^l)\right)\abs{z-z'};$$
\item for each $(t,x,x',y,z) \in [0,T] \times \R^d \times \R^d \times \R \times \R^{1\times d}$,
$$\abs{ f(t,x,y,z)-f(t,x',y,z)} \leqslant \left(C+\frac{\beta}{2}(\abs{x}^{1/l}+\abs{x'}^{1/l})\right)\abs{x-x'},$$
$$\abs{g(x)-g(x')} \leqslant \left(C+\frac{\alpha}{2}(\abs{x}^{1/l}+\abs{x'}^{1/l})\right)\abs{x-x'};$$
\item $$\alpha+T\beta < \frac{1}{e^{1/l}2^{1-1/l}\gamma^{1/l} e^{((1+1/l)K_b+K_{f,y}) T} \abs{\sigma}_{\infty}^{1+1/l} T^{1/l}}.$$
\end{enumerate}

Sometimes we will also consider stronger assumption.
\paragraph{Assumption (B.2).}
Let $f: [0,T] \times \mathbb{R}^d \times \mathbb{R} \times \mathbb{R}^{1\times d} \rightarrow \mathbb{R}$ and $g:\mathbb{R}^d \rightarrow  \mathbb{R}$ be continuous functions and let us assume moreover that there exist five constants, $l \geqslant 1$, $0 \leqslant r < \frac{1}{l}$, $\alpha \geqslant 0$, $\beta \geqslant 0$, $\gamma \geqslant 0$ and $K_{f,y} \geqslant 0$ such that:
\begin{enumerate}
\item for each $(t,x,y,y',z) \in [0,T] \times \mathbb{R}^d \times \mathbb{R} \times \mathbb{R} \times \mathbb{R}^{1\times d}$,
$$ \abs{f(t,x,y,z)-f(t,x,y',z)} \leqslant K_{f,y} \abs{y-y'};$$
\item for each $(t,x,y,z,z') \in [0,T] \times \R^d \times \R \times \R^{1\times d} \times \R^{1\times d}$,
$$\abs{ f(t,x,y,z)-f(t,x,y,z')} \leqslant \left(C+\frac{\gamma}{2}(\abs{z}^l+\abs{z'}^l)\right)\abs{z-z'};$$
\item for each $(t,x,x',y,z) \in [0,T] \times \R^d \times \R^d \times \R \times \R^{1\times d}$,
$$\abs{ f(t,x,y,z)-f(t,x',y,z)} \leqslant \left(C+\frac{\beta}{2}(\abs{x}^{r}+\abs{x'}^{r})\right)\abs{x-x'},$$
$$\abs{g(x)-g(x')} \leqslant \left(C+\frac{\alpha}{2}(\abs{x}^{r}+\abs{x'}^{r})\right)\abs{x-x'}.$$
\end{enumerate}

\begin{rem}
 Assumption (B.2) implies assumption (B.1). Moreover, the quadratic case corresponds to $l=1$.
\end{rem}

\begin{prop}
\label{proposition r<1/l}
We assume that assumptions (F.1) and (B.2) hold. There exists a solution $(Y,Z)$ of the Markovian BSDE (\ref{EDSR}) in $\mathcal{S}^2\times \mathcal{M}^2$ such that,
$$\abs{Z_t}  \leqslant C(1+\abs{X_t}^r).$$
Moreover, this solution is unique amongst solutions $(Y,Z)$ such that
\begin{itemize}
 \item $Y \in \mathcal{S}^2$,
 \item there exists $\eta >0$ such that
$$\mathbb{E} \left[e^{(\frac{1}{2}+\eta)\frac{\gamma^2}{4}\int_0^T \abs{Z_s}^{2l} ds}\right] < +\infty.$$
\end{itemize}
\end{prop}


\paragraph*{Proof of the proposition}
First of all, let us remark that if we have a solution $(Y,Z)$ such that 
$$\abs{Z_t} \leqslant C(1+\abs{X_t}^r), \quad \quad \forall t \in [0,T],$$
then, for all $c>0$,
\begin{equation}
\label{moment expo pour Z}
 \mathbb{E} \left[e^{c\int_0^T \abs{Z_s}^{2l} ds}\right] \leqslant \mathbb{E} \left[Ce^{C\sup_{0 \leqslant t \leqslant T} \abs{X_t}^{2lr} }\right]  < +\infty
\end{equation}
because $2lr<2$ (see e.g. part 5 in \cite{Briand-Hu-08}). Now, let us start by the uniqueness result. We consider two solutions
$(Y^1,Z^1)$ and $(Y^2,Z^2)$ such that $Y^1,Y^2 \in \mathcal{S}^2$, $\abs{Z^1} \leqslant C(1+\abs{X}^r)$ and there exists $\eta>0$ such that 
$$\mathbb{E} \left[e^{(\frac{1}{2}+\eta)\frac{\gamma^2}{4}\int_0^T \abs{Z_s^2}^{2l} ds}\right] < +\infty.$$
We define $\bar{Y}:=Y^1-Y^2$ and $\bar{Z} := Z^1-Z^2$. By considering the difference of the two BSDEs, the classical linearization method gives us
\begin{eqnarray*}
 \bar{Y}_t &=& \int_t^T \bar{Y}_s U_s+\bar{Z}_s V_s ds -\int_t^T \bar{Z}_s dW_s,
\end{eqnarray*}
that is to say
\begin{eqnarray}
\label{Ybar linearise}
\bar{Y}_t &=& -\int_t^T e^{\int_t^s U_udu}\bar{Z}_s (dW_s-V_sds),
\end{eqnarray}
where $(U,V)$ takes value in $\R \times \R^d$ and
$$\abs{U_s} \leqslant K_{f,y}, \quad \quad \quad \abs{V_s} \leqslant C+\frac{\gamma}{2}(\abs{Z^1_s}^l+\abs{Z^2_s}^l).$$
By applying Young's inequality, Hölder's inequality and (\ref{moment expo pour Z}), we have
\begin{eqnarray*}
 \mathbb{E} \left[e^{\frac{1}{2}\int_0^T \abs{V_s}^{2} ds}\right] &\leqslant& \E\left[e^{\frac{1}{2}\int_0^T (C+C\abs{Z_s^1}^{2l}+(1+\eta)\frac{\gamma^2}{4}\abs{Z_s^2}^{2l}) ds}\right]\\
 &\leqslant& C\E\left[e^{C\int_0^T \abs{Z_s^1}^{2l}ds}e^{\frac{1+\eta}{2}\frac{\gamma^2}{4}\int_0^T \abs{Z_s^2}^{2l} ds}\right]\\
 &\leqslant& C\E\left[e^{Cp\int_0^T \abs{Z_s^1}^{2l}ds}\right]^{1/p}\E\left[e^{(\frac{1}{2}+\eta)\frac{\gamma^2}{4}\int_0^T \abs{Z_s^2}^{2l} ds}\right]^{1/q}\\
 &<& +\infty,
\end{eqnarray*}
with $q=(1/2+\eta)(1/2+\eta/2)^{-1}$. This estimate shows us that Novikov's condition is fulfilled and so we are able to use Girsanov's Theorem in (\ref{Ybar linearise}) that gives us directly that $\bar{Y}=0$. Then it is standard to show that $\bar{Z}=0$. Finally we obtain the uniqueness result.

Now, let us show the existence result. Firstly we will approximate our Markovian BSDE by another one. Let $(Y^M,Z^M)$ the solution of the BSDE
\begin{equation}
\label{EDSR approchee}
Y^M_t = g_M(X_T)+\int_t^T f_M(s,X_s,Y^M_s,Z_s^M)ds-\int_t^T Z_s^M dW_s,
\end{equation}
with $g_M=g \circ \rho_M$ and $f_M=f(.,\rho_M(.),.,.)$ where $\rho_M$ is a smooth modification of the projection on the centered euclidean ball of radius $M$ such that $\abs{\rho_M}\leqslant M$, $\abs{\nabla \rho_M} \leqslant 1$ and $\rho_M(x)=x$ when $\abs{x}\leqslant M-1$. It is now easy to see that $g_M$ and $f_M$ are Lipschitz functions with respect to $x$. Theorem 3.1 in \cite{Richou-11} gives us that $Z^M$ is bounded by a constant $A_0$ that depends on $M$ in the quadratic case. In fact this result stays true in our more general framework. More precisely we have this proposition that we will show in the appendix.

\begin{prop}
\label{prop Z borne}$ $
We assume that (F.1) holds. We also assume that $f: [0,T] \times \mathbb{R}^d \times \mathbb{R} \times \mathbb{R}^{1\times d} \rightarrow \mathbb{R}$ and $g:\mathbb{R}^d \rightarrow  \mathbb{R}$ are continuous functions such that:
\begin{itemize}
 \item $g$ is $K_g$-Lipschitz,
 \item $f$ is $K_{f,x}$-Lipschitz with respect to $x$, $K_{f,y}$-Lipschitz with respect to $y$ and locally Lipschitz with respect to $z$: there exists an increasing function $\varphi : \R^+ \rightarrow \R^+$ such that for each $(t,x,y,z,z') \in [0,T] \times \R^d \times \R \times \R^{1\times d} \times \R^{1\times d}$,
$$\abs{ f(t,x,y,z)-f(t,x,y,z')} \leqslant C(1+\varphi(\abs{z})+\varphi(\abs{z'}))\abs{z-z'}.$$
Then, there exists a unique solution $(Y,Z)$ to the BSDE (\ref{EDSR}) in $\mathcal{S}^2 \times \mathcal{M}^2$ such that $Z$ is bounded. Moreover, we have
$$\abs{Z} \leqslant e^{(2K_b+K_{f,y})T}\abs{\sigma}_{\infty}(K_g+TK_{f,x}).$$
\end{itemize}
\end{prop}

Thanks to this lemma we know that there exists a unique solution $(Y^M,Z^M)$ to the BSDE (\ref{EDSR approchee}) (in the appropriate space) and $Z^M$ is bounded by a constant $A_0$ that depends on $M$. Moreover, $f_M$ is a Lipschitz function with respect to $z$ and BSDE (\ref{EDSR approchee}) is a classical Lipschitz BSDE. Now we will show the following lemma.

\begin{lem}
\label{lemme recurrence}
We have, 
 $$\abs{Z_t^M} \leqslant A_n+B_n\abs{X_t}^r,$$
 with $(A_n,B_n)_{n \in \N}$ defined by recursion: $B_0=0$, $A_0$ defined before,
$$B_{n+1}=C,$$
$$A_{n+1}=C(1+A_n^{rl}),$$
where $C$ is a constant that does not depend on $M$.
\end{lem}
\paragraph*{Proof of the lemma}
Let us prove the result by recursion. For $n=0$ we have already shown the result. Let us assume that the result is true for some $n \in \N^*$ and let us show that it stays true for $n+1$. Firstly we assume that for all $t \in [0,T]$, $b(t,.)$, $g$ and $f(t,.,.,.)$ are sufficiently differentiable. Then $X$ and $(Y^M,Z^M)$ are differentiable with respect to $x$ (see e.g. \cite{ElKaroui-Peng-Quenez-97}), we have
\begin{eqnarray*}
 \nabla Y_t^M &=& \nabla g_M(X_T) \nabla X_T-\int_t^T \nabla Z_s^M dW_s\\
 &&+ \int_t^T \nabla_x f_M(s,X_s,Y_s^M,Z_s^M)\nabla X_s+\nabla_y f_M(s,X_s,Y_s^M,Z_s^M)\nabla Y^M_s+\nabla_z f_M(s,X_s,Y_s^M,Z_s^M)\nabla Z_s^M ds,
\end{eqnarray*}
and $Z^M_t=\nabla Y^M_t (\nabla X_t)^{-1}\sigma(t) \textrm{ a.s.}$. Since $\abs{Z_s^M} \leqslant A_0$, we have
$$\abs{\nabla_z f_M(s,X_s,Y^M_s,Z_s^M)} \leqslant C(1+\abs{Z^M_s}^l) \leqslant C$$
and so we are allowed to apply Girsanov's Theorem: $\tilde{W}_t:=W_t-\int_0^t \nabla_z f_M(s,X_s,Y_s^M,Z_s^M) ds$ is a Brownian motion under a probability $\mathbb{Q}^M$. We obtain
\begin{eqnarray*}
\nabla Y_t^M &=&\mathbb{E}^{\mathbb{Q}^M}_t \Bigg[ e^{\int_t^T \nabla_y f_M(u,X_u,Y_u^M,Z_u^M) du}\nabla g_M(X_T)\nabla X_T\\
&&\left.+\int_t^T e^{\int_t^s \nabla_y f_M(u,X_u,Y_u^M,Z_u^M) du}\nabla_x f_M(s,X_s,Y_s^M,Z_s^M) \nabla X_s ds  \right], 
\end{eqnarray*}
and finally
\begin{equation}
\label{inegalite Z}
\abs{Z_t^M} \leqslant  C +e^{(K_b+K_{f,y})(T-t)}\abs{\sigma}_{\infty}\mathbb{E}^{\mathbb{Q}^M}_t \left[ \alpha\abs{X_T}^r +\beta \int_t^T \abs{X_s}^rds \right]
\end{equation}
because $\nabla X_s(\nabla X_t)^{-1}$ is bounded by $e^{K_b (T-t)}$. Let us come back to the SDE: we have
$$X_s=X_t+\int_t^s b(u,X_u)du + \int_t^s \sigma(u) d\tilde{W}_u + \int_t^s \sigma(u) \nabla_z f (u,X_u,Y_u^M,Z_u^M)du,$$
\begin{equation}
\label{point cle b sous lineaire}
\abs{X_s} \leqslant \abs{X_t} + C+\int_t^s K_b\abs{X_u}du + \abs{\int_t^s \sigma(u) d\tilde{W}_u}+ \abs{\sigma}_{\infty}\gamma\int_t^s \abs{A_n+B_n\abs{X_u}^r}^ldu,
\end{equation}
$$ \mathbb{E}^{\mathbb{Q}^M}_t \left[ \abs{X_s}  \right] \leqslant \abs{X_t} +C +K_b\int_t^s \mathbb{E}^{\mathbb{Q}^M}_t \left[ \abs{X_u}  \right]du + CA_n^l+ C\mathbb{E}^{\mathbb{Q}^M}_t \left[ \int_t^s B_n^l\abs{X_u}^{rl} du \right],$$
thanks to the recursion assumption. Young's inequality gives us
$$B_n^l\abs{X_u}^{rl} \leqslant \frac{B_n^{lp}}{p}+\frac{\abs{X_u}^{rlq}}{q}$$
with $1/p+1/q=1$ and $rlq=1$ (let us recall that $rl<1$). Thus, we obtain
$$ \mathbb{E}^{\mathbb{Q}^M}_t \left[ \abs{X_s}\right] \leqslant \abs{X_t} +C +K_b\int_t^s \mathbb{E}^{\mathbb{Q}^M}_t \left[ \abs{X_u}  \right]du + CA_n^l+ CB_n^{lp}+C\mathbb{E}^{\mathbb{Q}^M}_t \left[ \int_t^s \abs{X_u} du  \right].$$
Gronwall's Lemma gives us
$$\mathbb{E}^{\mathbb{Q}^M}_t \left[ \abs{X_s}  \right] \leqslant C+CA_n^l+CB_n^{lp}+C\abs{X_t},$$
and so
$$\mathbb{E}^{\mathbb{Q}^M}_t \left[ \abs{X_s}^r \right] \leqslant \left(\mathbb{E}^{\mathbb{Q}^M}_t \left[ \abs{X_s}  \right]\right)^r \leqslant C+CA_n^{rl}+CB_n^{rlp}+C\abs{X_t}^r$$
because $r<1$. By introducing this inequality into (\ref{inegalite Z}) we obtain
$$\abs{Z_t^M}  \leqslant C+CA_n^{rl}+CB_n^{rlp}+C\abs{X_t}^r.$$
Finally, we set
$$B_{n+1}=C,$$
and, thanks to the recursion assumption, we can take
$$A_{n+1}=C(1+ A_n^{rl}).$$
When $b$, $g$ and $f$ are not differentiable, we can prove the Lemma by a standard approximation and stability results for Lipschitz BSDEs.
\cqfd

Since $lr<1$, the recursion function that defines the sequence $(A_n)_{n \geqslant 0}$ is a contractor function and so $A_n \rightarrow A_{\infty}$ when $ n \rightarrow +\infty$, with $A_{\infty}$ that does not depend on $M$. Finally, 
\begin{equation}
\label{estime sur Z}
  \abs{Z_t^M}  \leqslant A_{\infty}+C\abs{X_t}^r.
\end{equation}

Now, we want to come back to the initial BSDE (\ref{EDSR}). We will show that $(Y^n, Z^n)_{n \in \mathbb{N}}$ is a Cauchy sequence in the space $\mathcal{S}^2 \times \mathcal{M}^2$. We have, thanks to the classical linearization method,
\begin{eqnarray*}
Y^{p+q}_t-Y^p_t &=& g_{p+q}(X_T)-g_p(X_T)+\int_t^T f_{p+q}(s,X_s,Y^{p+q}_s,Z^{p+q}_s)-f_{p}(s,X_s,Y^{p+q}_s,Z^{p+q}_s)ds\\
&&+\int_t^T(Y^{p+q}_s-Y^p_s)U^{p,q}_s+(Z^{p+q}_s-Z^p_s)V^{p,q}_sds-\int_t^T Z^{p+q}_s-Z^p_s dW_s,
\end{eqnarray*}
that is to say
\begin{eqnarray*}
 Y^{p+q}_t-Y^p_t &=& e^{\int_t^T U^{p,q}_udu}\left[g_{p+q}(X_T)-g_p(X_T)\right]\\
 &&+\int_t^Te^{\int_t^s U^{p,q}_udu}\left[ f_{p+q}(s,X_s,Y^{p+q}_s,Z^{p+q}_s)-f_{p}(s,X_s,Y^{p+q}_s,Z^{p+q}_s)\right]ds\\
 &&-\int_t^Te^{\int_t^s U^{p,q}_udu} (Z^{p+q}_s-Z^p_s)(dW_s-V^{p,q}_sds),
\end{eqnarray*}
with $p, q \in \mathbb{N}$, $\abs{U^{p,q}_s}\leqslant K_{f,y}$ and $\abs{V_s^{p,q}} \leqslant \frac{\gamma}{2}(1+\abs{Z^p_s}^l+\abs{Z^{p+q}_s}^l)$.  
Thanks to (\ref{estime sur Z}), Novikov's condition is fulfilled and so we are able to apply Girsanov's Theorem:
\begin{eqnarray*}
 \abs{Y^{p+q}_t-Y^p_t} &\leqslant& \abs{\mathbb{E}^{\mathbb{Q}^{p,q}}_t \left[ e^{\int_t^T U^{p,q}_udu}\left[g_{p+q}(X_T)-g_p(X_T)\right]  \right]}\\
 &&+\abs{\mathbb{E}^{\mathbb{Q}^{p,q}}_t \left[\int_t^T e^{\int_t^s U^{p,q}_udu}\left[ f_{p+q}(s,X_s,Y^{p+q}_s,Z^{p+q}_s)-f_{p}(s,X_s,Y^{p+q}_s,Z^{p+q}_s)\right]ds \right]}\\
 &\leqslant& e^{K_{f,y}T}\mathbb{E}^{\mathbb{Q}^{p,q}}_t \left[ \abs{g_{p+q}(X_T)-g_p(X_T)}  \right]\\
 &&+e^{K_{f,y}T}\int_t^T\mathbb{E}^{\mathbb{Q}^{p,q}}_t \left[ \abs{f_{p+q}(s,X_s,Y^{p+q}_s,Z^{p+q}_s)-f_{p}(s,X_s,Y^{p+q}_s,Z^{p+q}_s)}  \right]ds\\
 &\leqslant& C\mathbb{E}^{\mathbb{Q}^{p,q}}_t \left[ (1+\abs{X_T}^{r+1})\mathbbm{1}_{\abs{X_T}>p-1}  \right]+C\int_t^T\mathbb{E}^{\mathbb{Q}^{p,q}}_t \left[ (1+\abs{X_s}^{r+1})\mathbbm{1}_{\abs{X_s}>p-1} \right]ds\\
&\leqslant& C\mathbb{E}^{\mathbb{Q}^{p,q}}_t \left[ 1+\abs{X_T}^{2r+2}  \right]^{1/2} \mathbb{E}^{\mathbb{Q}^{p,q}}_t \left[\mathbbm{1}_{\abs{X_T}>p-1} \right]^{1/2} \\
&& +C \int_t^T \mathbb{E}^{\mathbb{Q}^{p,q}}_t \left[ 1+\abs{X_s}^{2r+2} \right]^{1/2} \mathbb{E}^{\mathbb{Q}^{p,q}}_t \left[\mathbbm{1}_{\abs{X_s}>p-1} \right]^{1/2}ds\\ 
 &\leqslant& C\mathbb{E}^{\mathbb{Q}^{p,q}}_t \left[ 1+\abs{X_T}^{2r+2} \right]^{1/2} \frac{\mathbb{E}^{\mathbb{Q}^{p,q}}_t \left[\abs{X_T} \right]^{1/2}}{p^{1/2}}\\
&& +C \int_t^T \mathbb{E}^{\mathbb{Q}^{p,q}}_t \left[ 1+\abs{X_s}^{2r+2}  \right]^{1/2} \frac{\mathbb{E}^{\mathbb{Q}^{p,q}}_t \left[\abs{X_s} \right]^{1/2}}{p^{1/2}}ds.
\end{eqnarray*}
By the same calculus than in the proof of Lemma (\ref{lemme recurrence}) using the a priori estimate on $Z^{p,q}$ and $Z^p$, we are able to show that
$$\mathbb{E}^{\mathbb{Q}^{p,q}}_t \left[\abs{X_s}^a \right] \leqslant C(1+\abs{X_t}^a)$$
for all $a \geqslant 1$ with a constant $C$ that depends on $a$ but does not depend on $p$ and $q$. Finally
$$\mathbb{E}\left[ \sup_{0 \leqslant t \leqslant T} \abs{Y^{p+q}_t-Y^p_t}^2 \right] \leqslant \frac{C(1+ \mathbb{E}\left[\sup_{0 \leqslant t \leqslant T} \abs{X_t}^{r+3/2}\right])}{p^{1/2}} \leqslant \frac{C}{p^{1/2}} \xrightarrow{ p \rightarrow +\infty} 0.$$
By applying Itô's formula to the process $\abs{Y^{p+q}-Y^p}^2$ and using the same calculus, it is rather standard to show that $(Z^n)_{n \in \mathbb{N}}$ is a Cauchy sequence in $\mathcal{M}^2$. Finally, it is easy to check that $(Y^n,Z^n) \xrightarrow{n \rightarrow +\infty} (Y,Z)$ the solution of the initial BSDE (\ref{EDSR}) and our estimate (\ref{estime sur Z}) on $Z^n$ stays true for $Z$:
$$\abs{Z_t} \leqslant A_{\infty}+C\abs{X_t}^r.$$ 
\cqfd

\begin{thm}
\label{theoreme r=1surl}
We assume that assumptions (F.1) and (B.1) hold. There exists a solution $(Y,Z)$ of the Markovian BSDE such that
$$\abs{Z_t} \leqslant C+e^{(K_b(1+1/l)+K_{f,y})(T-t)}(\alpha+\beta T)\abs{\sigma}_{\infty}e^{1/l}\abs{X_t}^{1/l}, \quad \forall t \in [0,T].$$
Moreover, this solution is unique amongst solutions $(Y,Z)$ such that 
\begin{itemize}
 \item $Y \in \mathcal{S}^2$,
 \item there exists $\eta >0$ such that
$$\mathbb{E} \left[e^{(2+\eta)\frac{\gamma^2}{4}\int_0^T \abs{Z_s}^{2l} ds}\right] < +\infty.$$
\end{itemize}
\end{thm}


\paragraph*{Proof of the theorem}
We will mimic the proof of Proposition (\ref{proposition r<1/l}). Let us start by the uniqueness result. We consider two solutions $(Y^1,Z^1)$ and $(Y^2,Z^2)$ such that $Y^1, Y^2 \in \mathcal{S}^2$ and 
$$\E \left[ e^{(2+\eta)\frac{\gamma^2}{4}\int_0^T \abs{Z_s^1}^{2l}ds}\right]+\E \left[e^{(2+\eta)\frac{\gamma^2}{4}\int_0^T \abs{Z_s^2}^{2l}ds}\right]<+\infty.$$
As in the proof of Proposition (\ref{proposition r<1/l}) we consider the BSDE satisfied by processes $\bar{Y}$ and $\bar{Z}$ and we introduce the two processes $U$ and $V$. Now we just have to show that Novikov's condition stays fulfilled: by applying Young's inequality and Hölder's inequality we have
\begin{eqnarray*}
  \E\left[e^{\frac{1}{2}\int_0^T \abs{V_s}^{2} ds}\right] &\leqslant& \mathbb{E} \left[e^{\frac{1}{2}\int_0^T (C+(2+\eta)\frac{\gamma^2}{4}\abs{Z_s^1}^{2l}+(2+\eta)\frac{\gamma^2}{4}\abs{Z_s^2}^{2l}) ds}\right]\\
 &\leqslant& C\E\left[e^{\frac{2+\eta}{2}\frac{\gamma^2}{4}\int_0^T \abs{Z_s^1}^{2l}ds}e^{\frac{2+\eta}{2}\frac{\gamma^2}{4}\int_0^T \abs{Z_s^2}^{2l} ds}\right]\\
 &\leqslant& C\E\left[e^{(2+\eta)\frac{\gamma^2}{4}\int_0^T \abs{Z_s^1}^{2l}ds}\right]^{1/2}\E\left[e^{(2+\eta)\frac{\gamma^2}{4}\int_0^T \abs{Z_s^2}^{2l} ds}\right]^{1/2}\\
 &<& +\infty.
\end{eqnarray*}
The remaining of the uniqueness proof is unchanged. Now, let us show the existence result. We will consider again the solution $(Y^M,Z^M)$ of the BSDE (\ref{EDSR approchee}). We have already remark that $Z^M$ is bounded by a constant $A_0$ that depends on $M$. Now we want to obtain an estimate on $Z^M$ that does not depend on $M$ by showing the following lemma.
\begin{lem}
\label{lemme recurrence2}
We have, 
 $$\abs{Z_t^M} \leqslant A_n(t)+B_n(t)\abs{X_t}^{1/l},$$
 with $(A_n,B_n)_{n \in \N}$ defined by recursion: $B_0=0$, $A_0$ defined before,
$$B_{n+1}(t)=\abs{\sigma}_{\infty}(\alpha+\beta T)e^{(K_b(1+1/l)+K_{f,y})(T-t)}e^{2^{l-1}\abs{\sigma}_{\infty}\gamma T B_n^l(t)/l},$$
$$A_{n+1}(t)=\abs{\sigma}_{\infty}(\alpha+\beta T)e^{(K_b(1+1/l)+K_{f,y})(T-t)}e^{2^{l-1}\abs{\sigma}_{\infty}\gamma T B_n^l(t)/l}\left(C+2^{1-1/l}\abs{\sigma}^{1/l}_{\infty}\gamma^{1/l} T^{1/l}A_n(t)\right),$$
where $C$ is a constant that does not depend on $M$ and $t$.
\end{lem}
\paragraph*{Proof of the lemma}
Let us prove the result by recursion. For $n=0$ we have already shown the result. Let us assume that the result is true for some $n \in \N^*$ and let us show that it stays true for $n+1$. Firstly we assume that for all $t \in [0,T]$, $b(t,.)$, $g$ and $f(t,.,.,.)$ are sufficiently differentiable. Then, by the same calculus than in the proof of Lemma \ref{lemme recurrence}, inequalities (\ref{inegalite Z}) and (\ref{point cle b sous lineaire}) stay true and we easily obtain
$$ \mathbb{E}^{\mathbb{Q}^M}_t \left[ \abs{X_s}  \right] \leqslant \abs{X_t} +C +K_b\int_t^s \mathbb{E}^{\mathbb{Q}^M}_t \left[ \abs{X_u}  \right]du + 2^{l-1}\abs{\sigma}_{\infty}\gamma TA_n^l(t)+ 2^{l-1}\abs{\sigma}_{\infty}\gamma B_n^l(t)\mathbb{E}^{\mathbb{Q}^M}_t \left[ \int_t^s \abs{X_u} du \right],$$
because $t\mapsto B_n(t)$ and $t\mapsto A_n(t)$ are not increasing functions. Gronwall's Lemma gives us
$$\mathbb{E}^{\mathbb{Q}^M}_t \left[ \abs{X_s}  \right] \leqslant \left(C+2^{l-1}\abs{\sigma}_{\infty}\gamma TA_n^l(t)+\abs{X_t}\right)e^{\left(K_b+2^{l-1}\abs{\sigma}_{\infty}\gamma B_n^l(t)\right)(T-t)},$$
and so
$$\mathbb{E}^{\mathbb{Q}^M}_t \left[ \abs{X_s}^{1/l} \right] \leqslant \left(\mathbb{E}^{\mathbb{Q}^M}_t \left[ \abs{X_s}  \right]\right)^{1/l} \leqslant \left(C+2^{1-1/l}\abs{\sigma}^{1/l}_{\infty}\gamma^{1/l} T^{1/l}A_n(t)+\abs{X_t}^{1/l}\right)e^{\frac{1}{l}\left(K_b+2^{l-1}\abs{\sigma}_{\infty}\gamma B_n^l(t)\right)(T-t)}$$
because $1/l<1$. By introducing this inequality into (\ref{inegalite Z}) we obtain
$$\abs{Z_t^M}  \leqslant C+\abs{\sigma}_{\infty}(\alpha+\beta T)e^{\left(K_b(1+1/l)+K_{f,y}+\frac{2^{l-1}\abs{\sigma}_{\infty}\gamma }{l}B_n^l(t)\right)(T-t)}\left(C+2^{1-1/l}\abs{\sigma}^{1/l}_{\infty}\gamma^{1/l} T^{1/l}A_n(t)+\abs{X_t}^{1/l}\right).$$
Finally, we can take
$$B_{n+1}(t)=\abs{\sigma}_{\infty}(\alpha+\beta T)e^{(K_b(1+1/l)+K_{f,y})(T-t)}e^{2^{l-1}\abs{\sigma}_{\infty}\gamma T B_n^l(t)/l},$$
and, 
$$A_{n+1}(t)=\abs{\sigma}_{\infty}(\alpha+\beta T)e^{(K_b(1+1/l)+K_{f,y})(T-t)}e^{2^{l-1}\abs{\sigma}_{\infty}\gamma T B_n^l(t)/l}\left(C+2^{1-1/l}\abs{\sigma}^{1/l}_{\infty}\gamma^{1/l} T^{1/l}A_n(t)\right).$$
When $b$, $g$ and $f$ are not differentiable, we can prove the lemma by a standard approximation and stability results for Lipschitz BSDEs.
\cqfd
Now we want to study the behavior of the sequence $(B_n(t))_{n \in \mathbb{N}}$. Let us denote
$$C_1(t):=\abs{\sigma}_{\infty}(\alpha+\beta T)e^{(K_b(1+1/l)+K_{f,y})(T-t)} \quad \textrm{ and } \quad C_2:= 2^{l-1}\abs{\sigma}_{\infty}\gamma T.$$
Then, we have $B_{n+1}(t)=C_1(t)e^{\frac{C_2}{l}B_n^l(t)}$. It is easy to see that this sequence has a finite limit $B_{\infty}(t) \in \R^+$ if and only if $\set{x\geqslant 0 | x=C_1(t)e^{\frac{C_2}{l}x^l}}$ is not empty and in this case we have
$$B_{\infty}(t)=\inf \set{x\geqslant 0 | x=C_1(t)e^{\frac{C_2}{l}x^l}}.$$
Moreover, the set $\set{x\geqslant 0 | x=C_1(t)e^{\frac{C_2}{l}x^l}}$ has only one, two or three elements. If $\set{x\geqslant 0 | x=C_1(t)e^{\frac{C_2}{l}x^l}}$ has only one element denoted by $\tilde{x}$, necessarily we have
$$
\left\{ \begin{array}{l} \tilde{x}=C_1(t)e^{\frac{C_2}{l}\tilde{x}^l}\\
1=C_1(t)C_2\tilde{x}^{l-1}e^{\frac{C_2}{l}\tilde{x}^l},
\end{array}
\right.
$$
that gives us
$$
\left\{ \begin{array}{l} \tilde{x}=C_1(t)e^{1/l}\\
C_1(t)=\frac{1}{(eC_2)^{1/l}}.
\end{array}
\right.
$$
Since we assume assumption (B.1).4, we have
$$C_1(t) \leqslant C_1(0) < \frac{1}{(eC_2)^{1/l}},$$
so we conclude that the sequence $(B_n(t))_{n \in \mathbb{N}}$ has a finite limit $B_{\infty}(t) <C_1(t)e^{1/l}$ that does not depend on $M$. Now, let us see what happens to the sequence 
$(A_n(t))_{n \in \mathbb{N}}$. We can remark that
$$A_{n+1}(t)=B_{n+1}(t)\left(C+C_2^{1/l}A_n(t)\right).$$
Since $B_n(t) \rightarrow B_{\infty}(t)<C_1(t)e^{1/l}$ and $C_1(t)e^{1/l}C_2^{1/l}<1$ due to assumption (B.1).4, a classical result for sequences defined by recursion gives us that $(A_n(t))_{n \in \mathbb{N}}$ has a finite limit $A_{\infty}(t)= \frac{CB_{\infty}(t)}{1-B_{\infty}(t)C_2^{1/l}}$ that does not depend on $M$ because it does not depend on the initial value $A_0(t)$.
Finally, we have
\begin{eqnarray}
 \nonumber \abs{Z^M_t} &\leqslant&  A_{\infty}(t)+B_{\infty}(t)\abs{X}^{1/l} \leqslant \frac{CB_{\infty}(0)}{1-B_{\infty}(0)C_2^{1/l}}+ C_1(t)e^{1/l}\abs{X}^{1/l}\\
\label{estimee ZM independant M} &\leqslant& C+\abs{\sigma}_{\infty}(\alpha+\beta T)e^{(K_b(1+1/l)+K_{f,y})(T-t)}e^{1/l}\abs{X}^{1/l}. 
\end{eqnarray}

To conclude we have to show that $(Y^n,Z^n)_{n \in \N}$ is a Cauchy sequence in the space $\mathcal{S}^2 \times \mathcal{M}^2$. As in the proof of Proposition~\ref{proposition r<1/l} we consider the BSDE satisfied by processes $Y^{p+q}-Y^p$ and $Z^{p+q}-Z^p$ and we introduce two processes $U^{p,q}$ and $V^{p,q}$. Now we just have to show that Novikov's condition stays fulfilled: by applying Young's inequality and Hölder's inequality we have
\begin{eqnarray}
\label{calcul condition novikov 1}
 \mathbb{E} \left[e^{\frac{1}{2}\int_0^T \abs{V_s^{p,q}}^{2} ds}\right] &\leqslant& \E\left[e^{\frac{1}{2}\int_0^T (C+(2+\eta)\frac{\gamma^2}{4}\abs{Z_s^{p+q}}^{2l}+(2+\eta)\frac{\gamma^2}{4}\abs{Z_s^{p}}^{2l}) ds}\right]\\
&\leqslant& C\E\left[e^{\frac{2+\eta}{2}\frac{\gamma^2}{4}\int_0^T \abs{Z_s^{p+q}}^{2l}ds}e^{\frac{2+\eta}{2}\frac{\gamma^2}{4}\int_0^T \abs{Z_s^p}^{2l} ds}\right]\\
\label{calcul condition novikov 3}
 &\leqslant& C\E\left[e^{(2+\eta)\frac{\gamma^2}{4}\int_0^T \abs{Z_s^{p+q}}^{2l}ds}\right]^{1/2}\E\left[e^{(2+\eta)\frac{\gamma^2}{4}\int_0^T \abs{Z_s^p}^{2l} ds}\right]^{1/2}.
\end{eqnarray}
Thanks to estimate (\ref{estimee ZM independant M}) we have, for all $M \in \N$,
\begin{eqnarray}
\nonumber &&\E\left[e^{(2+\eta)\frac{\gamma^2}{4}\int_0^T \abs{Z_s^{M}}^{2l}ds}\right]\\
\label{inegalite expo ZM Novikov} &\leqslant& C\E \left[ \exp\left((1+\eta)\frac{\gamma^2}{2}e^{2l(K_b(1+1/l)+K_{f,y})T}(\alpha+\beta T)^{2l}\abs{\sigma}_{\infty}^{2l}e^{2}T \sup_{0\leqslant t \leqslant T} e^{-2K_bt}\abs{X_t}^2\right) \right].
\end{eqnarray}
Since we have
$$\sup_{0\leqslant t \leqslant T} \left(e^{-2K_bt}\abs{X_t}^2\right) \leqslant \sup_{0\leqslant t \leqslant T} \left(e^{-2K_bt}\sup_{0\leqslant s \leqslant t}\abs{X_s}^2\right),$$
a slight modification of the proof of Lemma 4.1 in \cite{Delbaen-Hu-Richou-09} gives us that the right term in inequality (\ref{inegalite expo ZM Novikov}) is finite when
$$(1+\eta)\frac{\gamma^2}{2}e^{2l(K_b(1+1/l)+K_{f,y})T}(\alpha+\beta T)^{2l}\abs{\sigma}_{\infty}^{2l}e^{2}T < \frac{1}{2\abs{\sigma}_{\infty}^2 T}$$ 
which is true when $\eta$ is small enough because we assume in assumption (B.1) that
$$\alpha+T\beta < \frac{1}{e^{1/l}2^{1-1/l}\gamma^{1/l} e^{((1+1/l)K_b+K_{f,y}) T} \abs{\sigma}_{\infty}^{1+1/l} T^{1/l}},$$
that is to say,
$$\frac{\gamma^2}{2}e^{2l(K_b(1+1/l)+K_{f,y})T}(\alpha+\beta T)^{2l}\abs{\sigma}_{\infty}^{2l}e^{2}T < \frac{1}{2\abs{\sigma}_{\infty}^2 T}\frac{1}{2^{2(l-1)}}<\frac{1}{2\abs{\sigma}_{\infty}^2 T}.$$ 
Finally Novikov's condition is fulfilled and we are allowed to use Girsanov's theorem. The remaining of the existence proof is unchanged. To conclude, we have to show that when we have a solution that verifies
$$\abs{Z_t} \leqslant C+e^{(K_b(1+1/l)+K_{f,y})(T-t)}(\alpha+\beta T)\abs{\sigma}_{\infty}e^{1/l}\abs{X_t}^{1/l}, \quad \forall t \in [0,T],$$
then, there exists $\eta>0$ such that 
$$\mathbb{E} \left[e^{(2+\eta)\frac{\gamma^2}{4}\int_0^T \abs{Z_s}^{2l} ds}\right] < +\infty.$$
This claim is easy to prove because we have
\begin{eqnarray*}
\E\left[e^{(2+\eta)\frac{\gamma^2}{4}\int_0^T \abs{Z_s}^{2l}ds}\right] \leqslant C\E \left[ \exp\left((1+\eta)\frac{\gamma^2}{2}e^{2l(K_b(1+1/l)+K_{f,y})T}(\alpha+\beta T)^{2l}\abs{\sigma}_{\infty}^{2l}e^{2}T \sup_{0\leqslant t \leqslant T} e^{-2K_bt}\abs{X_t}^2\right) \right],
\end{eqnarray*}
and we have already shown that the right term in previous inequality is finite for $\eta$ small enough.
\cqfd

\begin{rem}
 In (B.1), our assumption 
$$\alpha+T\beta < \frac{1}{e\gamma^{1/l} e^{((1+l^{-1})K_b+K_{f,y}) T} \abs{\sigma}_{\infty}^{1+1/l} T^{1/l}}$$
is more restrictive than the one we can find in the article \cite{Delbaen-Hu-Richou-09} for the quadratic case (i.e. $l=1$), that is to says
$$\alpha+T\beta < \frac{1}{\gamma^{1/l} e^{((1+l^{-1})K_b+K_{f,y}) T} \abs{\sigma}_{\infty}^{1+1/l} T^{1/l}}.$$
In this case, the assumption in \cite{Delbaen-Hu-Richou-09} is more or less optimal because we need it to obtain a sufficient exponential moment for the terminal condition and the random part of the generator. Let us also remark that in \cite{Delbaen-Hu-Richou-09} assumptions are more general because they are about the growth of $f$ and $g$ instead of the growth of derivatives of $f$ and $g$.
\end{rem}

\begin{rem}
 With the same machinery it is possible to treat a little more general framework than the one of assumption (B.1): indeed it is possible to replace points 2 and 3 with
\begin{enumerate}
\setcounter{enumi}{1} 
\item for each $(t,x,y,z,z') \in [0,T] \times \R^d \times \R \times \R^{1\times d} \times \R^{1\times d}$,
$$\abs{ f(t,x,y,z)-f(t,x,y,z')} \leqslant \left(C+\gamma'\abs{x}+\frac{\gamma}{2}(\abs{z}^l+\abs{z'}^l)\right)\abs{z-z'};$$
\item for each $(t,x,x',y,z) \in [0,T] \times \R^d \times \R^d \times \R \times \R^{1\times d}$,
$$\abs{ f(t,x,y,z)-f(t,x',y,z)} \leqslant \left(C+\beta'\abs{z}+\frac{\beta}{2}(\abs{x}^{1/l}+\abs{x'}^{1/l})\right)\abs{x-x'},$$
$$\abs{g(x)-g(x')} \leqslant \left(C+\frac{\alpha}{2}(\abs{x}^{1/l}+\abs{x'}^{1/l})\right)\abs{x-x'};$$
\end{enumerate}
and the point 4 with an Ad hoc assumption. We decided to do not deal with this little more general setting because the proof is already technical and we do not want to complicate it unnecessarily.
\end{rem}

\section{Some results when $\sigma$ is random}
The main restriction in the previous part is about the function $\sigma$ that is assumed to be deterministic. In this section we will give some partial results when the SDE is given by
\begin{equation}
 \label{EDS sigma aleatoire}
X_t=x+\int_0^t b(s,X_s)ds+\int_0^t \sigma(s,X_s)dW_s.
\end{equation}
We will consider classical assumptions on this SDE.
\paragraph{Assumption (F.2).}
Let $b : [0,T] \times \mathbb{R}^d \rightarrow \mathbb{R}^d$ and $\sigma : [0,T]\times \R^d \rightarrow \mathbb{R}^{d \times d}$ be continuous  functions and let us assume that there exist $K_b \geqslant 0$, $K_{\sigma}\geqslant 0$ and $M_{\sigma}\geqslant 0$ such that:
\begin{enumerate}
 \item $\forall t \in [0,T]$, $\abs{b(t,0)} \leqslant C$,
 \item $\forall t \in [0,T]$, $\forall (x,x') \in \mathbb{R}^d \times \mathbb{R}^d$, $\abs{b(t,x)-b(t,x')} \leqslant K_b \abs{x-x'}$,
 \item $\forall t \in [0,T]$, $\forall x \in \mathbb{R}^d$, $\abs{\sigma(t,x)} \leqslant M_{\sigma}$,
 \item $\forall t \in [0,T]$, $\forall (x,x') \in \mathbb{R}^d \times \mathbb{R}^d$, $\abs{\sigma(t,x)-\sigma(t,x')} \leqslant K_{\sigma} \abs{x-x'}$.
\end{enumerate}
Before giving our first result, let us point out why we are not able to use the same machinery than in our first part. When $\sigma$ is deterministic, the starting point is Proposition~\ref{prop Z borne} where we show that $Z$ is bounded under good assumptions. To prove this result we deeply use the fact that $\nabla X$ is bounded. Now, when $\sigma$ is not deterministic, $\nabla X$ is not necessarily bounded and so Proposition 2.3 does not necessarily remain. Finally, the first question to answer is: does the process $Z$ remains bounded when $g$ and $f$ are Lipschitz with respect to $x$?

\subsection{Boundedness of $Z$ when $T$ is small enough}
In this part we will give a partial answer to the previous question.
\begin{prop}
 We assume that (F.2) holds. We also assume that $f: [0,T] \times \mathbb{R}^d \times \mathbb{R} \times \mathbb{R}^{1 \times d} \rightarrow \mathbb{R}$ and $g:\mathbb{R}^d \rightarrow  \mathbb{R}$ are continuous functions such that:
\begin{itemize}
 \item $g$ is $K_g$-Lipschitz,
 \item $f$ is $K_{f,x}$-Lipschitz with respect to $x$, $K_{f,y}$-Lipschitz with respect to $y$ and locally Lipschitz with respect to $z$: there exists an increasing continuous function $\varphi : \R^+ \rightarrow \R^+$ such that for each $(t,x,y,z,z') \in [0,T] \times \R^d \times \R \times \R^{1\times d} \times \R^{1\times d}$,
$$\abs{ f(t,x,y,z)-f(t,x,y,z')} \leqslant (K_{f,z}+\varphi(\abs{z})+\varphi(\abs{z'}))\abs{z-z'}.$$
Then, for $T$ small enough, there exists a unique solution $(Y,Z)$ to the BSDE (\ref{EDSR}) in $\mathcal{S}^2 \times \mathcal{M}^2$ such that $Z$ is bounded.
\end{itemize}
\end{prop}

\paragraph*{Proof of the proposition}
Once again, we will use a classical truncation argument (see e.g. the proof of Theorem 4.1 in \cite{Delbaen-Hu-Bao-09}). Our truncation function $\rho_M$ is a smooth modification of the projection on the centered euclidean ball of radius $M$ such that $\abs{\rho_M}\leqslant M$, $\abs{\nabla \rho_M} \leqslant 1$ and $\rho_M(x)=x$ when $\abs{x}\leqslant M-1$. We denote $(Y^M,Z^M)$ the solution of the BSDE
$$Y_t^M=g(X_T)+\int_t^T f_M(s,X_s,Y_s^M,Z_s^M)ds-\int_t^T Z_s^MdW_s,$$
where $f_M:=f(.,.,.,\rho_M(.))$. Now, this BSDE is also Lipschitz with respect to $z$. Firstly we assume that for all $t \in [0,T]$, $b(t,.)$, $g$ and $f(t,.,.,.)$ are differentiable. Then $X$ and $(Y,Z)$ are differentiable with respect to $x$, we have
\begin{eqnarray*}
 \nabla Y_t^M &=& \nabla g(X_T) \nabla X_T-\int_t^T \nabla Z_s^M dW_s\\
 &&+ \int_t^T \nabla_x f_M(s,X_s,Y_s^M,Z_s^M)\nabla X_s+\nabla_y f_M(s,X_s,Y_s^M,Z_s^M)\nabla Y^M_s+\nabla_z f_M(s,X_s,Y_s^M,Z_s^M)\nabla Z_s^M ds,
\end{eqnarray*}
and $Z^M_t=\nabla Y^M_t (\nabla X_t)^{-1}\sigma(t,X_t) \textrm{ a.s.}$. Since $\nabla_z f_M$ is bounded by $K_{f,z}+2\varphi(M)$, we are allowed to apply Girsanov's Theorem: $\tilde{W}_t:=W_t-\int_0^t \nabla_z f_M(s,X_s,Y_s^M,Z_s^M) ds$ is a Brownian motion under the probability $\mathbb{Q}^M$. We obtain
\begin{eqnarray*}
\nabla Y_t^M &=&\mathbb{E}^{\mathbb{Q}^M}_t \Bigg[ e^{\int_t^T \nabla_y f_M(u,X_u,Y_u^M,Z_u^M) du}\nabla g_M(X_T)\nabla X_T\\
&&\left.+\int_t^T e^{\int_t^s \nabla_y f_M(u,X_u,Y_u^M,Z_u^M) du}\nabla_x f_M(s,X_s,Y_s^M,Z_s^M) \nabla X_s ds \right], 
\end{eqnarray*}
and finally
$$\abs{Z_t^M} \leqslant e^{K_{f,y}T}M_{\sigma}\left( K_g \E^{\Q^M}_t \left[\abs{U_T}\right] + K_{f,x}\int_t^T\E^{\Q^M}_t \left[ \abs{U_s}\right]ds\right),$$
with $(U_s)_{t \leqslant s\leqslant T}$ the solution to the SDE
\begin{eqnarray*}
 U_s &=& Id+\int_t^s \nabla b(X_u)U_udu+\int_t^s \sum_{i=1} \nabla \sigma^i (X_u)U_u dW_u^i\\
 & = & Id+\int_t^s \nabla b(X_u)U_udu+\int_t^s \sum_{i=1} \nabla \sigma^i (X_u)U_u (d\tilde{W}_u^i +(\nabla_z f_M)^i(u,X_u,Y_u^M,Z_u^M) du)
\end{eqnarray*}
where the superscript $^i$ denotes the $i$-th column. A classical estimate on $\E^{\Q^M}_t \left[\abs{U_s}\right]$ gives us 
$$\E^{\Q^M}_t \left[\abs{U_s}\right] \leqslant Ce^{K_{\sigma}^2T+\left(K_b+K_{\sigma}(K_{f,z}+2\varphi(M))\right)^2T^2}.$$
For $T$ small enough, the function 
$$x \mapsto 1+Ce^{K_{f,y}T}M_{\sigma}( K_g + K_{f,x}T)e^{K_{\sigma}^2T+\left(K_b+K_{\sigma}(K_{f,z}+2\varphi(x))\right)^2T^2}$$
has, at least, one positive fixed point. Let us denotes $\bar{M}$ the lowest positive fixed point of this function. Then $Z^{\bar{M}} \leqslant \bar{M}-1$ and $(Y^{\bar{M}},Z^{\bar{M}})$ is a solution to the initial BSDE. Uniqueness follows from the uniqueness result for Lipschitz BSDEs.
 \cqfd 

\begin{rem}
$ $
 \begin{itemize}
  \item In general, it is not possible to stick local solutions to obtain a solution $(Y,Z)$ with $Z$ bounded for all $T$.
 \item The biggest $T$ that allows the existence of a fixed point for the function
 $$x \mapsto 1+Ce^{K_{f,y}T}M_{\sigma}( K_g + K_{f,x}T)e^{K_{\sigma}^2T+\left(K_b+K_{\sigma}(K_{f,z}+2\varphi(x))\right)^2T^2},$$
strongly depends on $K_g$. So, it is not possible to treat the case of $g$ and $f$ locally Lipschitz with respect to $x$ by using the same machinery than in the previous part.
 \end{itemize}
\end{rem}

\subsection{A simple example}
In this part we will see that when we consider a simple quadratic BSDE with an explicit solution, the process $Z$ remains bounded when $g$ and $f$ are Lipschitz with respect to $x$. More precisely, we will consider the following quadratic BSDE:
\begin{equation}
\label{exemple simple EDSR quadratique}
Y_t=g(X_T)+\int_t^T \frac{\abs{Z_s}^2}{2}ds-\int_t^T Z_sdW_s.
\end{equation}
\begin{prop}
 Let us assume that (F.2) holds and that $g$ is a $K_g$-Lipschitz function. Then there exists a unique solution $(Y,Z)$ to the BSDE (\ref{exemple simple EDSR quadratique}) in $\mathcal{S}^2 \times \mathcal{M}^2$ such that the process $Z$ is bounded.
\end{prop}
\paragraph{Proof of the proposition}
It is well known that (\ref{exemple simple EDSR quadratique}) can be explicitly solved with an exponential transform, also called Cole-Hopf transform in PDE theory. More precisely, we have
$$Y_t=\log \E_t \left[e^{g(X_T)}\right], \quad Z_t=e^{-Y_t}\tilde{Z}_t,$$
where $\tilde{Z}$ is given by the martingale representation theorem applied to the martingale $\left(\E_t\left[e^{g(X_T)}\right]\right)_{0 \leqslant t \leqslant T}$. $Y$ is well defined because $g$ is Lipschitz and for all $C>0$
$$\E \left[e^{C\abs{X_T}}\right]<+\infty,$$
since $\sigma$ is bounded (see e.g. part 5 in \cite{Briand-Hu-08}). The uniqueness is standard. As in previous proofs, we assume in a first time that $g$, $b$ and $\sigma$ are differentiable with respect to $x$. Then we have
\begin{eqnarray*}
 \abs{Z_t} &=& \abs{\frac{\E_t\left[\nabla g(X_T)\nabla X_T (\nabla X_t)^{-1}e^{g(X_T)}\right]}{\E_t\left[e^{g(X_T)}\right]}\sigma(t,X_t)}\\
&\leqslant& C\frac{\E_t\left[\abs{\nabla X_T (\nabla X_t)^{-1}}e^{g(X_T)}\right]}{\E_t\left[e^{g(X_T)}\right]}\\
&\leqslant& C\E_t\left[\abs{\nabla X_T (\nabla X_t)^{-1}}^2\right]^{1/2}\frac{\E_t\left[e^{2g(X_T)}\right]^{1/2}}{\E_t\left[e^{g(X_T)}\right]}\\
&\leqslant& C\frac{\E_t\left[e^{2g(X_T)}\right]^{1/2}}{\E_t\left[e^{g(X_T)}\right]}
\end{eqnarray*}
because $\nabla g$, $\sigma$ are bounded and $(\nabla X_s (\nabla X_t)^{-1})_{t \leqslant s \leqslant T}$ solve the SDE
$$U_s=Id+\int_t^s \nabla b(u,X_u)U_u du+\sum_{i=1}^d \int_t^s \nabla \sigma^i (u,X_u)U_u dW_u^i,$$
so $\E_t\left[\abs{\nabla X_T (\nabla X_t)^{-1}}^2\right]$ is bounded. Let us denote $\bar{X}$ the solution of the ordinary differential equation (with a random initial condition)
$$\bar{X}_s=X_t+\int_t^s b(u,\bar{X}_u)du, \quad t\leqslant s \leqslant T.$$
Since $g$ is Lipschitz and $\bar{X}_T$ is $\mathcal{F}_t$-measurable, we obtain
\begin{equation}
\label{majoration Z exemple}
\abs{Z_t}\leqslant C\frac{e^{g(\bar{X}_T)}\E_t\left[e^{2C\abs{X_T-\bar{X}_T}}\right]^{1/2}}{e^{g(\bar{X}_T)}\E_t\left[e^{-C\abs{X_T-\bar{X}_T}}\right]} \leqslant C\frac{\E_t\left[e^{2C\abs{X_T-\bar{X}_T}}\right]^{1/2}}{\E_t\left[e^{-C\abs{X_T-\bar{X}_T}}\right]}.
\end{equation}
Let us estimate $\E_t\left[e^{2C\abs{X_T-\bar{X}_T}}\right]$. We have
$$\abs{X_s-\bar{X}_s} = 0+K_b\int_t^s \abs{X_u-\bar{X}_u}du +\sup_{t \leqslant r \leqslant T}\abs{\int_t^r \sigma(u,X_u)dW_u},$$
and we deduce from Gronwall's lemma the inequality
$$ \sup_{t \leqslant s \leqslant T}\abs{X_s-\bar{X}_s} \leqslant C\sup_{t \leqslant s \leqslant T}\abs{\int_t^s \sigma(u,X_u)dW_u}.$$
It follows from the Dambis-Dubins-Schwarz representation theorem that, for $\lambda \geqslant 0$, 
\begin{eqnarray*}
 \mathbb{E}\left[ \exp \left(\lambda\sup_{t \leqslant s \leqslant T}\abs{\int_t^s \sigma(u,X_u)dW_u}\right)\bigg| \mathcal{F}_t, X_t=x_0 \right] &=& \mathbb{E}_t\left[ \exp \left(\lambda\sup_{t \leqslant s \leqslant T}\abs{\int_t^s \sigma(u,X^{t,x_0}_u)dW_u}\right)\right]\\
& \leqslant & \E\left[ \sup_{0 \leqslant s \leqslant \norm{\sigma}_{\infty}^2T} e^{ \lambda\abs{W_s}} \right]<+\infty
\end{eqnarray*}
where $(X_s^{t,x_0})_{t \leqslant s \leqslant T}$ stands for the solution to the SDE (\ref{EDS sigma aleatoire}) that starts from $x_0$ at time $t$. We remark that the right term in the last inequality does not depend on $x_0$ so 
$$\mathbb{E}_t\left[ \exp \left(C\sup_{t \leqslant s \leqslant T}\abs{\int_t^s \sigma(u,X_u)dW_u}\right) \right]$$
is upper bounded. By the same type of argument we have that
$$\mathbb{E}_t\left[ \exp \left(-C\sup_{t \leqslant s \leqslant T}\abs{\int_t^s \sigma(u,X_u)dW_u}\right) \right]$$
is lower bounded by a strictly positive constant. Finally (\ref{majoration Z exemple}) becomes
$$\abs{Z_t} \leqslant C\frac{\E_t\left[e^{2C\sup_{t \leqslant s \leqslant T}\abs{\int_t^s \sigma(u,X_u)dW_u}}\right]^{1/2}}{\E_t\left[e^{-C\sup_{t \leqslant s \leqslant T}\abs{\int_t^s \sigma(u,X_u)dW_u}}\right]}\leqslant C,$$
and so $Z$ is bounded. Finally, when $g$, $b$ and $\sigma$ are not differentiable, we can prove the result by a standard approximation. \cqfd
\begin{rem}
 Thanks to this estimate on $Z$, it is possible to use the same machinery than in the previous section to show estimates on $Z$ when $g$ and $f$ are locally Lipschitz with respect to $x$. This simple example is a good argument to postulate that Theorem \ref{theoreme r=1surl} or Proposition \ref{proposition r<1/l} could stay true when we replace (F.1) by (F.2), at least in the quadratic case. 
\end{rem}

\subsection{The case of bounded terminal conditions}
In this part we will restrict our study to the quadratic case and we will assume that the terminal condition and the generator are bounded with respect to $x$. In this case we are able to obtain estimates on $Z$ thanks to the additional tool of Bounded Mean Oscillation martingales (BMO martingales for short). We refer the reader to \cite{Kazamaki-94} for the theory of BMO martingales and we just recall the properties that we will use in the sequel. Let $\Phi_t=\int_0^t \phi_s dW_s$, for $t \in [0,T]$, be a real square integrable martingale with respect to the Brownian filtration. Then $\Phi$ is a BMO martingale if
$$\norm{\Phi}_{BMO}=\sup_{\tau \in [0,T]} \E\left[\langle \Phi \rangle_T-\langle \Phi \rangle_{\tau}|\mathcal{F}_{\tau}\right]^{1/2}=\sup_{\tau \in [0,T]} \E\left[\int_{\tau}^T \phi_s^2ds\bigg|\mathcal{F}_{\tau}\right]^{1/2} < +\infty,$$
where the supremum is taken over all stopping times in $[0,T]$ and $\langle \Phi \rangle$ denotes the quadratic variation of $\Phi$. In our case, the very important feature of BMO martingales is the following lemma:
\begin{lem}
\label{propriétés martingales BMOs}
 Let $\Phi$ be a BMO martingale. Then we have:
\begin{enumerate}
 \item The stochastic exponential 
$$\mathcal{E}(\Phi)_t=\mathcal{E}_t=\exp \left(\int_0^t \phi_s dW_s-\frac{1}{2}\int_0^t \abs{\phi_s}^2ds \right), \quad 0 \leqslant t \leqslant T,$$
is a uniformly integrable martingale.
 \item Thanks to the reverse H\"older inequality, there exists $p>1$ such that $\mathcal{E}_T \in L^p$. The maximal $p$ with this property can be expressed in terms of the BMO norm of $\Phi$.
\end{enumerate}
\end{lem}
We will work under following assumptions on coefficients of SDE (\ref{EDSR approchee cas quadratique sigma aleatoire}) and BSDE (\ref{EDSR}).
\paragraph{Assumption (F.3).}
Let $b : [0,T] \times \mathbb{R}^d \rightarrow \mathbb{R}^d$ and $\sigma : [0,T]\times \R^d \rightarrow \mathbb{R}^{d \times d}$ be continuous  functions and let us assume that there exist $K_b \geqslant 0$, $K_{\sigma}\geqslant 0$, $M_{\sigma}\geqslant 0$ and $\kappa \in[0,1]$ such that:
\begin{enumerate}
 \item $\forall t \in [0,T]$, $\abs{b(t,0)} \leqslant C$,
 \item $\forall t \in [0,T]$, $\forall (x,x') \in \mathbb{R}^d \times \mathbb{R}^d$, $\abs{b(t,x)-b(t,x')} \leqslant K_b \abs{x-x'}$,
 \item $\forall t \in [0,T]$, $\forall x \in \mathbb{R}^d$, $\abs{\sigma(t,x)} \leqslant M_{\sigma}(1+\abs{x}^{\kappa})$,
 \item $\forall t \in [0,T]$, $\forall (x,x') \in \mathbb{R}^d \times \mathbb{R}^d$, $\abs{\sigma(t,x)-\sigma(t,x')} \leqslant K_{\sigma} \abs{x-x'}$.
\end{enumerate}
\paragraph{Assumption (B.3).}
Let $f: [0,T] \times \mathbb{R}^d \times \mathbb{R} \times \mathbb{R}^{1\times d} \rightarrow \mathbb{R}$ and $g:\mathbb{R}^d \rightarrow  \mathbb{R}$ be continuous functions and let us assume moreover that there exist seven constants, $r \in \R^+$, $\alpha \geqslant 0$, $\beta \geqslant 0$, $\gamma \geqslant 0$, $K_{f,y} \geqslant 0$, $M_f \geqslant 0$ and $M_g \geqslant 0$ such that:
\begin{enumerate}
\item for each $(t,x,y,y',z) \in [0,T] \times \mathbb{R}^d \times \mathbb{R} \times \mathbb{R} \times \mathbb{R}^{1\times d}$,
$$ \abs{f(t,x,y,z)-f(t,x,y',z)} \leqslant K_{f,y} \abs{y-y'};$$
\item for each $(t,x,y,z,z') \in [0,T] \times \R^d \times \R \times \R^{1\times d} \times \R^{1\times d}$,
$$\abs{ f(t,x,y,z)-f(t,x,y,z')} \leqslant \left(C+\frac{\gamma}{2}(\abs{z}+\abs{z'})\right)\abs{z-z'};$$
\item for each $(t,x,x',y,z) \in [0,T] \times \R^d \times \R^d \times \R \times \R^{1\times d}$,
$$\abs{ f(t,x,y,z)-f(t,x',y,z)} \leqslant \left(C+\frac{\beta}{2}(\abs{x}^{r}+\abs{x'}^{r})\right)\abs{x-x'},$$
$$\abs{g(x)-g(x')} \leqslant \left(C+\frac{\alpha}{2}(\abs{x}^{r}+\abs{x'}^{r})\right)\abs{x-x'};$$
\item for each $(t,x,y,z) \in [0,T] \times \R^d  \times \R \times \R^{1\times d}$,
$$\abs{ f(t,x,y,z)} \leqslant M_f(1+\abs{y}+\abs{z}^2),$$
$$\abs{g(x)} \leqslant M_g.$$
\end{enumerate}
\begin{thm}
\label{theoreme g borne}
We assume that assumptions (F.3) and (B.3) hold. There exists a solution $(Y,Z)$ of the Markovian BSDE in $\mathcal{S}^2\times \mathcal{M}^2$ and this solution is unique amongst solutions $(Y,Z) \in \mathcal{S}^2\times \mathcal{M}^2$ such that $Y$ is bounded. Moreover we have
$$\abs{Z_t} \leqslant C(1+\abs{X_t}^{r+\kappa}), \quad \forall t \in [0,T],$$
and 
$$\norm{\int_0^. Z_sdW_s}_{BMO} \leqslant C,$$
where the last constant $C$ depends only on $M_g$, $M_f$ and $K_{f,y}$.
\end{thm}
\paragraph{Proof of the theorem}
For the existence and uniqueness result we refer the reader to \cite{Kobylanski-00,Lepeltier-SanMartin-98}. The estimate for the BMO norm of $Z$ is shown in \cite{Briand-Confortola-08,Ankirchner-Imkeller-Reis-07}. It just remains to prove the estimate on $Z$. As in previous proofs, we firstly assume that $f$, $g$, $b$ and $\sigma$ are differentiable with respect to $x$. Then, according to \cite{Briand-Confortola-08,Ankirchner-Imkeller-Reis-07}, $X$ and $(Y,Z)$ are differentiable with respect to $x$, we have
\begin{eqnarray*}
 \nabla Y_t &=& \nabla g(X_T) \nabla X_T-\int_t^T \nabla Z_s dW_s\\
 &&+ \int_t^T \nabla_x f(s,X_s,Y_s,Z_s)\nabla X_s+\nabla_y f(s,X_s,Y_s,Z_s)\nabla Y_s+\nabla_z f(s,X_s,Y_s,Z_s)\nabla Z_s ds,
\end{eqnarray*}
and $Z_t=\nabla Y_t (\nabla X_t)^{-1}\sigma(t,X_t) \textrm{ a.s.}$. Since $\int_0^. Z_sdW_s$ is BMO and 
$$\abs{\nabla_z f(s,X_s,Y_s,Z_s)} \leqslant C(1+\abs{Z_s})$$
then $\int_0^. \nabla_z f(s,X_s,Y_s,Z_s)dW_s$ is BMO and we are allowed to apply Girsanov's Theorem thanks to Lemma \ref{propriétés martingales BMOs}: $\tilde{W}_t:=W_t-\int_0^t \nabla_z f(s,X_s,Y_s,Z_s) ds$ is a Brownian motion under the probability 
$$\mathbb{Q}=\mathcal{E}\left(\int_0^. \nabla_z f(s,X_s,Y_s,Z_s)dW_s\right)_T\P.$$ 
We obtain
\begin{eqnarray*}
\nabla Y_t &=&\mathbb{E}^{\mathbb{Q}}_t \Bigg[ e^{\int_t^T \nabla_y f(u,X_u,Y_u,Z_u) du}\nabla g(X_T)\nabla X_T\\
&&\left.+\int_t^T e^{\int_t^s \nabla_y f(u,X_u,Y_u,Z_u) du}\nabla_x f(s,X_s,Y_s,Z_s) \nabla X_s ds \right], 
\end{eqnarray*}
and then it comes
\begin{equation}
 \label{equation intermediaire g bornee}
\abs{Z_t} \leqslant  C\left(1 +\mathbb{E}^{\mathbb{Q}}_t \left[ \abs{X_T}^{2r} \right]^{1/2} + \int_t^T \mathbb{E}^{\mathbb{Q}}_t \left[\abs{X_s}^{2r} \right]^{1/2}ds\right)\mathbb{E}^{\mathbb{Q}}_t \left[ \sup_{t \leqslant s \leqslant T}\abs{\nabla X_s(\nabla X_t)^{-1}}^{2} \right]^{1/2}\left(1+\abs{X_t}^{\kappa}\right)
\end{equation}
by using assumptions (F.3), (B.3) and Cauchy-Schwarz's inequality. Let us denote 
$$\mathcal{E}_{t,T}:=\exp \left(\int_t^T \nabla_z f(s,X_s,Y_s,Z_s) dW_s-\frac{1}{2}\int_t^T \abs{\nabla_z f(s,X_s,Y_s,Z_s)}^2ds \right).$$
Thanks to Lemma \ref{propriétés martingales BMOs}, there exists $p>1$ (that does not depend on $t$) such that $\E_t [\mathcal{E}_{t,T}^p ]<+\infty$. But, by using Hölder's inequality and classical estimates on SDEs we have
\begin{equation*}
 \mathbb{E}^{\mathbb{Q}}_t \left[ \abs{X_s}^{2r} \right] \leqslant \E_t [\mathcal{E}_{t,T}^p ]^{1/p}\mathbb{E}_t \left[ \abs{X_s}^{2rq} \right]^{1/q}
\leqslant C(1+\abs{X_t}^{2r}),
\end{equation*}
and
\begin{equation*}
 \mathbb{E}^{\mathbb{Q}}_t \left[ \sup_{t \leqslant s \leqslant T}\abs{\nabla X_s(\nabla X_t)^{-1}}^{2} \right] \leqslant \E_t [\mathcal{E}_{t,T}^p ]^{1/p}\mathbb{E}_t \left[ \sup_{t \leqslant s \leqslant T}\abs{\nabla X_s(\nabla X_t)^{-1}}^{2q} \right]^{1/q} \leqslant C,
\end{equation*}
By putting th two last inequalities into (\ref{equation intermediaire g bornee}) we obtain the result. Finally, when $b$, $g$ and $f$ are not differentiable, we can prove the result by standard approximations and stability results for quadratic BSDEs (see e.g. \cite{Imkeller-dosReis-09}). \cqfd

\section{Application to quadratic and superquadratic PDEs}
In this section we give an application of previous results concerning BSDEs to semilinear PDEs which have a quadratic or superquadratic growth with respect to the gradient of the solution. We will restrict our study to deterministic functions $\sigma$. Let us consider the following semilinear PDE
\begin{equation}
\label{EDP}
\left\{
\begin{array}{l}
\partial_t u(t,x)+\mathcal{L}u(t,x)+f(t,x,u(t,x),^t\nabla u(t,x)\sigma(t))=0, \quad x \in \R^d, t \in [0,T],\\
u(T,.)=g,
\end{array}
\right.
\end{equation}
where $\mathcal{L}$ is the infinitesimal generator of the diffusion $X^{t,x}$ solution to the SDE
$$X_s^{t,x}=x+\int_t^s b(r,X_r^{t,x})dr+\int_t^s \sigma(r)dW_r, \quad t \leqslant s \leqslant T.$$
The nonlinear Feynman-Kac formula consists in proving that the function defined by the formula
\begin{equation}
\label{definition u}
\forall (t,x) \in [0,T] \times \R^d, \quad u(t,x):=Y_t^{t,x}
\end{equation}
where, for each $(t_0,x_0) \in [0,T] \times \R^d$, $(Y^{t_0,x_0},Z^{t_0,x_0})$ stands for the solution given by Theorem~\ref{theoreme r=1surl} to the following BSDE
$$Y_t^{t_0,x_0}=g(X_T^{t_0,x_0})+\int_t^T f(s,X_s^{t_0,x_0},Y_s^{t_0,x_0},Z_s^{t_0,x_0})ds-\int_t^T Z_s^{t_0,x_0}dW_s, \quad 0 \leqslant t \leqslant T,$$
is a solution, at least a viscosity solution, to the PDE~(\ref{EDP}).
Firstly, let us study the growth and the continuity of this function.
\begin{prop}
\label{u continue}
Let assumptions (F.1) and (B.1) hold. The function $u$ defined by (\ref{definition u}) has a polynomial growth and is a continuous function. More precisely we have, $\forall (t,t',x,x') \in [0,T]^2\times \R^d\times \R^d$,
$$\abs{u(t,x)} \leqslant C(1+\abs{x}^{1+1/l}),$$
$$\abs{u(t,x)-u(t',x')} \leqslant C(1+\abs{x}^{1/l}+\abs{x'}^{1/l})\abs{x-x'}+C(1+\abs{x}^{1+1/l}+\abs{x'}^{1+1/l})\abs{t-t'}^{1/2}.$$
\end{prop}
\paragraph{Proof of the proposition}
To show the first point, it is sufficient to prove the estimate
\begin{equation}
 \label{estimee sur moment sup Y}
\E \left[\sup_{t \leqslant s \leqslant T} \abs{Y_s^{t,x}}^2 \right] \leqslant C \left(1+\abs{x}^{2(1+1/l)}\right).
\end{equation}
By a very classical method we can easily show the estimate
\begin{eqnarray*}
 \E \left[\sup_{t \leqslant s \leqslant T} \abs{Y_s^{t,x}}^2 \right] &\leqslant& C\E\left[ \abs{g(X_T^{t,x})}^2+\int_t^T \abs{f(s,X_s^{t,x},0,Z_s^{t,x})}^2ds\right].
\end{eqnarray*}
Since $\abs{Z}^{t,x} \leqslant C(1+\abs{X^{t,x}}^{1/l})$, we obtain, by using the growth of $g$ and $f$ and classical estimates on SDEs,
\begin{eqnarray*}
 \E \left[\sup_{t \leqslant s \leqslant T} \abs{Y_s^{t,x}}^2 \right] &\leqslant& C\E\left[ 1+\sup_{t\leqslant s \leqslant T}\abs{X_s^{t,x}}^{2(1+1/l)}\right]\\
&\leqslant & C \left(1+\abs{x}^{2(1+1/l)}\right).
\end{eqnarray*}
Now, let us show the second part of the proposition. By a symmetry argument we are allowed to suppose that $t'\geqslant t$. Then
$$u(t,x)-u(t',x') =\E\left[ Y_t^{t,x}-Y_{t'}^{t,x}\right]+ \E\left[ Y_{t'}^{t,x}-Y_{t'}^{t',x'}\right].$$
Cauchy-Schwarz's inequality and growth assumptions on $f$ and $g$ give us
\begin{eqnarray*}
 \abs{\E \left[ Y_t^{t,x}-Y_{t'}^{t,x}\right] }^2 &=& \abs{\E \int_{t}^{t'} f(s,X_s^{t,x},Y_s^{t,x},Z_s^{t,x})ds }^2\\
&\leqslant & \abs{t-t'} \E\left[ \int_{t}^{t'} \abs{f(s,X_s^{t,x},Y_s^{t,x},Z_s^{t,x})}^2ds\right]\\
&\leqslant & C\abs{t-t'} \E\left[1+ \sup_{t \leqslant s \leqslant T} \left(\abs{X_s}^{2(1+1/l)} +\abs{Y_s}^2 +\abs{Z_s}^{2l}\right)\right].
\end{eqnarray*}
Thanks to a priori estimate on $Z$, a classical estimate on SDEs and (\ref{estimee sur moment sup Y}), we obtain
$$\abs{\E \left[ Y_t^{t,x}-Y_{t'}^{t,x}\right] }^2 \leqslant C\abs{t-t'}(1+\abs{x}^{2(1+1/l)}).$$
Now we will study the term $\E\left[ Y_{t'}^{t,x}-Y_{t'}^{t',x'}\right]$. We have, thanks to the classical linearization method,
\begin{eqnarray*}
 Y^{t,x}_{t'}-Y^{t',x'}_{t'} &=& e^{\int_{t'}^T U^{x,x'}_udu}\left[g(X_T^{t,x})-g(X_T^{t',x'})\right]\\
 &&+\int_{t'}^Te^{\int_{t'}^s U^{x,x'}_udu}\left[ f(s,X_s^{t,x},Y^{t,x}_s,Z^{t,x}_s)-f(s,X_s^{t',x'},Y^{t,x}_s,Z^{t,x}_s)\right]ds\\
 &&-\int_{t'}^Te^{\int_{t'}^s U^{x,x'}_udu} (Z^{t,x}_s-Z^{t',x'}_s)(dW_s-V^{x,x'}_sds),
\end{eqnarray*}
with $\abs{U^{x,x'}_s}\leqslant K_{f,y}$ and $\abs{V_s^{x,x'}} \leqslant \frac{\gamma}{2}(1+\abs{Z^{t,x}_s}^l+\abs{Z^{t',x'}_s}^l)$.
Since Novikov's condition is fulfilled, we are able to apply Girsanov's Theorem. We obtain, by using the fact that $f$ and $g$ are locally Lipschitz,
\begin{eqnarray*}
 \abs{Y^{t,x}_{t'}-Y^{t',x'}_{t'}} &\leqslant& C\E^{\Q^{x,x'}}_{t'}\left[\abs{g(X_T^{t,x})-g(X_T^{t',x'})}\right]\\
 &&+C\E^{\Q^{x,x'}}_{t'}\left[\int_{t'}^T\abs{f(s,X_s^{t,x},Y^{t,x}_s,Z^{t,x}_s)-f(s,X_s^{t',x'},Y^{t,x}_s,Z^{t,x}_s)}ds\right]\\
&\leqslant& C\E^{\Q^{x,x'}}_{t'}\left[\sup_{t'\leqslant s\leqslant T} \left(1+\abs{X_s^{t,x}}^{1/l}+\abs{X_s^{t',x'}}^{1/l}\right)\abs{X_s^{t,x}-X_s^{t',x'}}\right]\\
&\leqslant& C\left(1+\E^{\Q^{x,x'}}_{t'}\left[\sup_{t'\leqslant s\leqslant T} \abs{X_s^{t,x}}^{2/l}\right]^{1/2}+\E^{\Q^{x,x'}}_{t'}\left[\sup_{t'\leqslant s\leqslant T}\abs{X_s^{t',x'}}^{2/l}\right]^{1/2}\right)\\
&&\times\E^{\Q^{x,x'}}_{t'}\left[\sup_{t'\leqslant s\leqslant T}\abs{X_s^{t,x}-X_s^{t',x'}}^2\right]^{1/2}.
\end{eqnarray*}
Let us recall that $\abs{V_s^{x,x'}} \leqslant C\left(1+\abs{X_s^{t,x}}+\abs{X_s^{t',x'}}\right)$. Once again we are able to use classical methods on SDEs to obtain finally
$$\abs{Y^{t,x}_{t'}-Y^{t',x'}_{t'}} \leqslant C(1+\abs{X_{t'}^{t,x}}^{1/l}+\abs{x'}^{1/l})\left(\abs{X_{t'}^{t,x}-x'}+\abs{t-t'}^{1/2}(1+\abs{X_{t'}^{t,x}}+\abs{x'})\right).$$
Classical estimates on SDEs allow us to conclude:
$$\E\abs{Y^{t,x}_{t'}-Y^{t',x'}_{t'}} \leqslant C(1+\abs{x}^{1/l}+\abs{x'}^{1/l})\left(\abs{x-x'}+\abs{t-t'}^{1/2}(1+\abs{x}+\abs{x'})\right).$$
\cqfd

\begin{prop}
 Let assumptions (F.1) and (B.1) hold. The function $u$ defined by (\ref{definition u}) is a viscosity solution to the PDE (\ref{EDP}).
\end{prop}
Since we are able to use Girsanov's transformation in the BSDE, we have a comparison result. Moreover, Proposition~\ref{u continue} gives us that $u$ is a continuous function. So the proof of the proposition is totally standard: for example, it can be easily adapted from the proof of Theorem 4.2 in \cite{ElKaroui-Peng-Quenez-97}.

\section{Time approximation of quadratic and superquadratic Markovian BSDEs}
\subsection{approximation of the initial BSDE by a Lipschitz one}
In a first time, we will consider the deterministic case for the function $\sigma$ and we will approximate the solution $(Y,Z)$ of the BSDE (\ref{EDSR}) by $(Y^M,Z^M)$ the solution of the BSDE (\ref{EDSR approchee}). The aim of the following proposition is to study the approximation error given by:
\begin{equation}
\label{erreur 1}
e_1(M):= \E \left[\sup_{0 \leqslant t \leqslant T} \abs{Y_t-Y^M_t}^2 \right]+ \E \left[ \int_0^T \abs{Z_t-Z^M_t}^2dt \right].
\end{equation}
\begin{prop}
\label{proposition etude erreur 1}
If we assume that assumptions (F.1) and (B.1) hold, then there exists $\lambda >0$ such that
$$e_1(M) \leqslant Ce^{-\lambda M^2}.$$
\end{prop}

\paragraph{Proof of the proposition}
Let us define processes $\delta Y:=Y-Y^M$ and $\delta Z := Z-Z^M$. We have
$$\delta Y_t=g(X_T)-g_M(X_T)+\int_t^T f(s,X_s,Y_s,Z_s)- f_M(s,X_s,Y_s^M,Z_s^M)ds-\int_t^T \delta Z_s dW_s.$$
The classical linearization method gives us
\begin{equation}
\label{equation aux differences linearisee}
\delta Y_t = \delta g + \int_t^T \delta f_s+\delta Y_s U_s^M+\delta Z_s V_s^M ds-\int_t^T \delta Z_s dW_s,
\end{equation}
with 
$$\delta g := g(X_T)-g_M(X_T),\quad \delta f_s := f(s,X_s,Y_s,Z_s)-f_M(s,X_s,Y_s,Z_s),$$
$(U^M,V^M)$ with value in $\R \times \R^d$ and
$$\abs{U^M_s} \leqslant K_{f,y}, \quad \quad \quad \abs{V^M_s} \leqslant C+\frac{\gamma}{2}(\abs{Z_s}^l+\abs{Z^M_s}^l).$$
We can easily show that Novikov's condition is fulfilled for $V^M$ by doing the same calculus than for $V^{p,q}$ in the proof of Theorem~\ref{theoreme r=1surl} (inequalities (\ref{calcul condition novikov 1}) to (\ref{calcul condition novikov 3})). So, we are allowed to apply Girsanov's theorem: $\tilde{W}_t := W_t-\int_0^t V^M_sds$ is a Brownian motion under the probability $\mathbb{Q}^M$. Thus, by applying Cauchy-Schwarz's inequality and Markov's inequality we obtain
\begin{eqnarray}
\nonumber
\delta Y_t &=& \E^{\Q^M}_t \left[ e^{\int_t^T U_s^M ds} \delta g +\int_t^T e^{\int_t^s U_u^M du}\delta f_s ds\right],\\
\nonumber
 \abs{\delta Y_t} &\leqslant& C \E^{\Q^M}_t \left[(1+\abs{X_T}^{1+1/l})\mathbbm{1}_{\abs{X_T} \geqslant M} +\int_t^T (1+\abs{X_s}^{1+1/l})\mathbbm{1}_{\abs{X_s} \geqslant M}ds \right]\\
\nonumber
&\leqslant& C \left(1+\E^{\Q^M}_t \left[\abs{X_T}^{2(1+1/l)} \right]\right)^{1/2} \frac{\E^{\Q^M}_t \left[e^{2\lambda\abs{X_T}^{2}} \right]^{1/2}}{e^{\lambda M^{2}}}\\
\label{estimation delta Y}
&& + C \int_t^T\left(1+\E^{\Q^M}_t \left[\abs{X_s}^{2(1+1/l)} \right]\right)^{1/2} \frac{\E^{\Q^M}_t \left[e^{2\lambda \abs{X_s}^{2}} \right]^{1/2}}{e^{\lambda M^{2}}}ds.
\end{eqnarray}
Then we use the following lemma that we will prove in the appendix.
\begin{lem}
\label{lemme moment expo sous Q}
 We assume that assumptions (F.1) and (B.1) hold. We have 
\begin{itemize}
 \item $\forall a \in [1,+\infty[$, $\exists C>0$,
$$\E^{\Q^M}_t \left[ \sup_{t \leqslant s\leqslant T} \abs{X_s}^a \right] \leqslant C(1+\abs{X_t}^a),\quad \forall t \in [0,T],$$
 \item $\exists C>0$, $\exists \bar{\mu}>0$, $\forall \mu \in [0,\bar{\mu}[$,
$$\E^{\Q^M}_t \left[\sup_{t \leqslant s\leqslant T} e^{\mu\abs{X_s}^{2}} \right]\leqslant Ce^{C\mu\abs{X_t}^{2}}, \quad \forall t \in [0,T].$$
\end{itemize}

\end{lem}
Now (\ref{estimation delta Y}) becomes,  
$$\abs{\delta Y_t} \leqslant \frac{C(1+ \abs{X_t}^{1+1/l})}{e^{\lambda M^2}} e^{C\lambda\abs{X_t}^{2}}.$$
By using Cauchy-Schwarz's inequality, we obtain for all $p\geqslant 1$ and for all $0<\lambda <\bar{\lambda}$ with $\bar{\lambda}$ small enough,
\begin{eqnarray*}
 \E\left[ \sup_{0 \leqslant t \leqslant T} \abs{\delta Y_t}^p \right] &\leqslant& \frac{C}{e^{p\lambda M^{2}}}\E \left[ (1+\sup_{0 \leqslant t \leqslant T} \abs{X_t}^{p(1+1/l)}) e^{Cp\lambda\sup_{0 \leqslant t \leqslant T} \abs{X_t}^{2}} \right]\\
 &\leqslant & \frac{C}{e^{p\lambda M^{2}}}\left( 1+\E\left[ \sup_{0 \leqslant t \leqslant T} \abs{X_t}^{2p(1+1/l)} \right]^{1/2}\right) \E \left[ e^{Cp\lambda\sup_{0 \leqslant t \leqslant T} \abs{X_t}^2} \right]^{1/2}.
\end{eqnarray*}
Let us remark that $C$ depends on $\bar{\lambda}$ but does not depend on $\lambda$. By using classical results about SDEs (see e.g. the beginning of part~5 in \cite{Briand-Hu-08}) we have, for all $p\geqslant 1$,
$$E\left[ \sup_{0 \leqslant t \leqslant T} \abs{X_t}^{2p(1+1/l)} \right]^{1/2} <+\infty,$$
and, for $\lambda$ small enough,
$$\E \left[ e^{Cp\lambda\sup_{0 \leqslant t \leqslant T} \abs{X_t}^2} \right]^{1/2} < +\infty.$$
Finally we obtain that
\begin{equation}
 \label{estimation erreur Lp en Y}
\E\left[ \sup_{0 \leqslant t \leqslant T} \abs{\delta Y_t}^p \right] \leqslant  \frac{C}{e^{p\lambda M^{2}}}.
\end{equation}

To study the error on $Z$ we come back to (\ref{equation aux differences linearisee}) and we apply Itô's formula:
$$\abs{\delta Y_t}^2 +\int_t^T \abs{\delta Z_s}^2ds= \abs{\delta Y_T}^2+ \int_t^T 2\delta Y_s ( \delta f_s+\delta Y_s U_s^M+\delta Z_s V_s^M) ds-\int_t^T 2\delta Y_s\delta Z_s dW_s.$$
We obtain by applying Cauchy-Schwarz's inequality
\begin{eqnarray*}
 \E \left[\int_0^T\abs{\delta Z_t}^2dt \right] & \leqslant& \E \left[ \abs{\delta Y_T}^2 \right]+2\E\left[\int_0^T\delta Y_t\delta f_tdt \right] + 2\E\left[\int_0^T\abs{\delta Y_t}^2U_t^Mdt \right] + 2\E\left[\int_0^T\delta Y_t\delta Z_tV_t^Mdt \right]\\
&\leqslant& (1+2K_{f,y})\E\left[ \sup_{0 \leqslant t \leqslant T} \abs{\delta Y_t}^2 \right] +2 T^{1/2}\E\left[ \sup_{0 \leqslant t \leqslant T} \abs{\delta Y_t}^2 \right]^{1/2}\E\left[ \sup_{0 \leqslant t \leqslant T} \abs{\delta f_t}^2 \right]^{1/2}\\
& &+2\E \left[\int_0^T\abs{\delta Z_t}^2dt \right]^{1/2} \E\left[ \sup_{0 \leqslant t \leqslant T} \abs{\delta Y_t}^4 \right]^{1/4}\E\left[ \sup_{0 \leqslant t \leqslant T} \abs{\delta V_t^M}^4 \right]^{1/4}.
\end{eqnarray*}
Thanks to inequalities
$$\abs{\delta f_t} \leqslant C(1+\abs{X_t}^{1/l})$$
and
$$\abs{\delta V_t^M} \leqslant C+\frac{\gamma}{2}(\abs{Z_t}^l+\abs{Z^M_t}^l) \leqslant C(1+\abs{X_t}),$$
it is easy to see that 
$$\E\left[ \sup_{0 \leqslant t \leqslant T} \abs{\delta f_t}^2 \right] + \E\left[ \sup_{0 \leqslant t \leqslant T} \abs{\delta V_t^M}^4 \right] \leqslant C$$
with $C$ that does not depend on $M$. Then, by applying (\ref{estimation erreur Lp en Y}) and the inequality $2ab\leqslant \frac{a^2}{2}+2b^2$ we have
\begin{eqnarray*}
 \E \left[\int_0^T\abs{\delta Z_t}^2dt \right]&\leqslant& \frac{C}{e^{2\lambda M^{2}}} +\frac{C}{e^{\lambda M^{2}}}+2\E \left[\int_0^T\abs{\delta Z_t}^2dt \right]^{1/2} \frac{C}{e^{\lambda M^{2}}} \\
&\leqslant & \frac{C}{e^{\lambda M^{2}}} +\frac{1}{2}\E \left[\int_0^T\abs{\delta Z_t}^2dt \right].
\end{eqnarray*}
Finally we obtain
\begin{equation}
 \label{estimation erreur en Z}
\E \left[\int_0^T\abs{\delta Z_t}^2dt \right] \leqslant  \frac{C}{e^{p\lambda M^{2}}}.
\end{equation}
To conclude, (\ref{estimation erreur Lp en Y}) and (\ref{estimation erreur en Z}) give us the result.
\cqfd

Now we want to obtain the same type of estimate in the quadratic case when $\sigma$ is random. Since $\sigma$ is not necessarily bounded, $Z$ could be unbounded even if $g$ and $f$ are Lipschitz functions with respect to $x$. So, we will approximate the solution $(Y,Z)$ of the BSDE (\ref{EDSR}) by $(\bar{Y}^M,\bar{Z}^M)$ the solution of the BSDE 
\begin{equation}
\label{EDSR approchee cas quadratique sigma aleatoire} 
\bar{Y}^M_t = g(\rho_{M^{(r+\kappa)^{-1}}}(X_T))+\int_t^T f(s,\rho_{M^{(r+\kappa)^{-1}}}(X_s),\bar{Y}^M_s,\rho_M(\bar{Z}_s^M))ds-\int_t^T \bar{Z}_s^M dW_s
\end{equation}
where $\rho_M$ is a smooth modification of the projection on the centered euclidean ball of radius $M$ such that $\abs{\rho_M}\leqslant M$, $\abs{\nabla \rho_M} \leqslant 1$ and $\rho_M(x)=x$ when $\abs{x}\leqslant M-1$.
The aim of the following proposition is to study the approximation error given by:
\begin{equation}
\label{erreur 1 bar}
\bar{e}_1(M):= \E \left[\sup_{0 \leqslant t \leqslant T} \abs{Y_t-\bar{Y}^M_t}^2 \right]+ \E \left[ \int_0^T \abs{Z_t-\bar{Z}^M_t}^2dt \right].
\end{equation}
\begin{prop}
\label{proposition etude erreur 1 bar}
If we assume that assumptions (F.3), (B.3) hold and $2\kappa \leqslant 1-r$, then there exists $\lambda >0$ such that
$$\bar{e}_1(M) \leqslant Ce^{-\lambda M^2}.$$
If moreover $2\kappa < 1-r$, then there exist $\lambda >0$ and $\varepsilon>0$ such that
$$\bar{e}_1(M) \leqslant Ce^{-\lambda M^{2+\varepsilon}}.$$
\end{prop}

\paragraph{Proof of the proposition}
Thanks to BMO tool, we have a comparison theorem that gives us an estimate for $\bar{e}_1$. Indeed, we can apply Lemma 3.2 in \cite{Imkeller-dosReis-09}: there exists $q>1$ such that
\begin{eqnarray}
\label{equation comparaison erreur approximation cas quadratique}
\nonumber
 \bar{e}_1(M) &\leqslant& C\E\left[\abs{g(X_T)-g(\rho_{M^{(1-\kappa)^{-1}}}(X_T))}^{2q}\right]^{1/q}\\
&&+C\E\left[\left( \int_0^T \abs{f(s,X_s,Y_s,Z_s)-f(s,\rho_{M^{(1-\kappa)^{-1}}}(X_s),Y_s,\rho_{M}(Z_s))}ds\right)^{2q}\right]^{1/q}.
\end{eqnarray}
Assumptions (F.3), (B.3) and the estimate on $Z$ give us
\begin{eqnarray*}
\abs{f(s,X_s,Y_s,Z_s)-f(s,\rho_{M^{(r+\kappa)^{-1}}}(X_s),Y_s,\rho_{M}(Z_s))} &\leqslant& C(1+\abs{X_s}^r)\mathbbm{1}_{\abs{X_s}\geqslant M^{(r+\kappa)^{-1}}} + C(1+\abs{Z_s})\mathbbm{1}_{\abs{Z_s}\geqslant M}\\
&\leqslant& C(1+\abs{X_s}^r)\mathbbm{1}_{\abs{X_s}^{r+\kappa}\geqslant M}\\
&& + C(1+\abs{X_s}^{r+\kappa})\mathbbm{1}_{\abs{X_s}^{r+\kappa}\geqslant M/C-1}\\
&\leqslant& C(1+\sup_{0 \leqslant s \leqslant T}\abs{X_s}^{r+\kappa})\mathbbm{1}_{\sup_{0 \leqslant s \leqslant T}\abs{X_s}^{r+\kappa}\geqslant M/C-1}
\end{eqnarray*}
and
\begin{eqnarray*}
\abs{g(X_T)-g(\rho_{M^{(1-\kappa)^{-1}}}(X_T))} &\leqslant& (1+\abs{X_T}^r)\mathbbm{1}_{\abs{X_T}\geqslant M^{(r+\kappa)^{-1}}}\\
&\leqslant&C(1+\sup_{0 \leqslant s \leqslant T}\abs{X_s}^{r+\kappa})\mathbbm{1}_{\sup_{0 \leqslant s \leqslant T}\abs{X_s}^{r+\kappa}\geqslant M/C-1}.
\end{eqnarray*}
By using Hölder's inequality and the fact that, for all $p>0$,
$$\E\left[\sup_{0 \leqslant s \leqslant T}\abs{X_s}^p\right]<+\infty,$$
(\ref{equation comparaison erreur approximation cas quadratique}) becomes
\begin{equation}
\label{equation comparaison erreur approximation cas quadratique 2}
 \bar{e}_1(M) \leqslant C\P\left(\sup_{0 \leqslant s \leqslant T}\abs{X_s}^{r+\kappa}\geqslant M/C-1 \right)^{q'}
\end{equation}
with $q'>1$. To conclude, we will use the following lemma that will be proved in the appendix.
\begin{lem}
\label{moment expo X sigma non borne}
 We assume that (F3) holds. There exists $\lambda>0$ such that
\begin{equation*}
 \mathbb{E} \left[ \exp \left( \lambda \sup_{0 \leqslant t \leqslant T} \abs{X_t}^{2(1-\kappa)} \right) \right] < +\infty.
\end{equation*}
\end{lem}
Since we have assume that $r+\kappa \leqslant 1-\kappa$, Markov's inequality and previous lemma give us, for $M$ big enough,
\begin{eqnarray*}
 \P\left(\sup_{0 \leqslant s \leqslant T}\abs{X_s}^{r+\kappa}\geqslant M/C-1 \right) &\leqslant& \frac{e^{\lambda\sup_{0 \leqslant s \leqslant T}\abs{X_s}^{2(1-\kappa)}}}{e^{\lambda(M/C-1)^2}}\leqslant \frac{C}{e^{\tilde{\lambda}M^2}}.
\end{eqnarray*}
Then, the first part of the Lemma is obtained by putting this inequality in the estimate (\ref{equation comparaison erreur approximation cas quadratique 2}). When $r+\kappa < 1-\kappa$, we denote $\varepsilon := \frac{1-\kappa}{r+\kappa}-1>0$ and we obtain, for $M$ big enough,
\begin{eqnarray*}
 \P\left(\sup_{0 \leqslant s \leqslant T}\abs{X_s}^{r+\kappa}\geqslant M/C-1 \right) &\leqslant& \frac{e^{\lambda\sup_{0 \leqslant s \leqslant T}\abs{X_s}^{2(1-\kappa)}}}{e^{\lambda(M/C-1)^{2(1+\varepsilon)}}}\leqslant \frac{C}{e^{\tilde{\lambda}M^{2(1+\varepsilon)}}}.
\end{eqnarray*}
Finally, the last part of the Lemma is proved by putting this inequality in the estimate (\ref{equation comparaison erreur approximation cas quadratique 2}). \cqfd

\begin{rem}
\label{remarque 2k >1-r}
 When $2\kappa > 1-r$ it is possible to show with the same proof that there exists $\lambda>0$ such that
 $$\bar{e}_1(M) \leqslant \frac{C}{\exp\left(\lambda M^{2\frac{1-\kappa}{r+\kappa}}\right)}.$$
When $\kappa=1$, it is also possible to recover the result obtained by Imkeller and dos Reis in \cite{Imkeller-dosReis-09} (they assume in addition that $r=0$): for all $k \in \N$, there exists $C>0$ such that
 $$\bar{e}_1(M) \leqslant \frac{C}{M^k}.$$
Let us remark that our result is more precise than the one of \cite{Imkeller-dosReis-09} and our proof is more simple since we do not have to study the second order Malliavin differentiability of the BSDE.
\end{rem}

\subsection{time approximation of the BSDE}
In a second time, we will approximate our modified BSDE by a discrete time one. We denote the time step by $h=\frac{T}{n}$ and $(t_k=kh)_{0 \leqslant k\leqslant n}$ stand for the discretization times. One needs to approximate $X$ by a Markov chain $X^n$ which can be simulated. For example, we will consider the classical Euler scheme given by 
\begin{equation*}
 \left\{
\begin{array}{lcl}
X_0^n & = & x\\
X^n_{t_{k+1}} & = & X^n_{t_{k}}+hb(t_k,X^n_{t_k})+\sigma(t_k,X^n_{t_k})(W_{t_{k+1}}-W_{t_k}), \quad 0 \leqslant k \leqslant n-1.
\end{array}
\right.
\end{equation*}
We denote $(Y^{M,n},Z^{M,n})$ (resp. $(\bar{Y}^{M,n},\bar{Z}^{M,n})$) our time approximation of $(Y^M,Z^M)$ (resp. $(\bar{Y}^M,\bar{Z}^M)$). These couples are obtained by the classical explicit dynamic programming equation:

\begin{equation*}
 \left\{
\begin{array}{lcl}
Y_{t_n}^{M,n} & = & g_M(X_{t_n}^n),\\
Z_{t_{k}}^{M,n} & = & \frac{1}{h}\E_{t_k} \left[Y_{t_{k+1}}^{M,n}(W_{t_{k+1}}-W_{t_k})\right], \quad 0 \leqslant k \leqslant n,\\
Y_{t_k}^{M,n} &=& \E_{t_k} \left[Y_{t_{k+1}}^{M,n}+hf_M(t_k,X^n_{t_k},Y^{M,n}_{t_{k+1}},Z^{M,n}_{t_k})\right], \quad 0 \leqslant k \leqslant n,
\end{array}
\right.
\end{equation*}
and
\begin{equation*}
 \left\{
\begin{array}{lcl}
\bar{Y}_{t_n}^{M,n} & = & g(\rho_{M^{(r+\kappa)^{-1}}}(X_{t_n}^n)),\\
\bar{Z}_{t_{k}}^{M,n} & = & \frac{1}{h}\E_{t_k} \left[\bar{Y}_{t_{k+1}}^{M,n}(W_{t_{k+1}}-W_{t_k})\right], \quad 0 \leqslant k \leqslant n,\\
\bar{Y}_{t_k}^{M,n} &=& \E_{t_k} \left[\bar{Y}_{t_{k+1}}^{M,n}+hf(t_k,\rho_{M^{(r+\kappa)^{-1}}}(X^n_{t_k}),\bar{Y}^{M,n}_{t_{k+1}},\rho_M(\bar{Z}^{M,n}_{t_k}))\right], \quad 0 \leqslant k \leqslant n.
\end{array}
\right.
\end{equation*}
In a classical framework, there is already results about the speed of convergence of BSDE time approximation. Let us precise the classical result shown by \cite{Bouchard-Touzi-04,Zhang-04,Lemor-05}.
\begin{prop}
 \label{prop resultat classique discretisation temporelle EDSR}
Let us assume that assumption (F.1) or (F.3) holds. We also assume that
\begin{itemize}
 \item $g$ is $K_g$-Lipschitz,
 \item $f$ is $K_{f,x}$-Lipschitz with respect to $x$, Lipschitz with respect to $y$ and $K_{f,z}$-Lipschitz with respect to $z$.
\end{itemize}
We denote $(Y^n,Z^n)$ the time discretization of $(Y,Z)$ given by the classical explicit dynamic programming equation. The error discretization is given by
$$e(n):=\sup_{0\leqslant k\leqslant n} \E \left[ \abs{Y^n_{t_k}-Y_{t_k}}^2 \right] + \sum_{k=0}^{n-1} \E \left[ \int_{t_k}^{t_{k+1}} \abs{Z^n_{t_k}-Z_t}^2dt\right].$$
Then, there exists a constant $C$ that does not depend on $K_g$, $K_{f,x}$ and $K_{f,z}$ such that
$$e(n) \leqslant Ce^{CK_{f,z}^2}\left[1+K_g^2+K_{f,x}^2+\E\left[\int_0^T \abs{Z_t}^2dt \right]+\E \left[ \sup_{0 \leqslant t \leqslant T} \abs{Y_t}^2\right]\right]h.$$
\end{prop}
 This proposition will be proved in the appendix. Now, the aim of this section is to study errors of discretization $e(M,n)$ and $\bar{e}(M,n)$ given by
$$e(M,n):= \sup_{0\leqslant k\leqslant n} \E \left[ \abs{Y^{M,n}_{t_k}-Y_{t_k}}^2 \right] + \sum_{k=0}^{n-1} \E \left[ \int_{t_k}^{t_{k+1}} \abs{Z^{M,n}_{t_k}-Z_t}^2dt\right]$$
and
$$\bar{e}(M,n):= \sup_{0\leqslant k\leqslant n} \E \left[ \abs{\bar{Y}^{M,n}_{t_k}-Y_{t_k}}^2 \right] + \sum_{k=0}^{n-1} \E \left[ \int_{t_k}^{t_{k+1}} \abs{\bar{Z}^{M,n}_{t_k}-Z_t}^2dt\right].$$

\begin{thm}
\label{theorem vitesse convergence sigma deterministe}
We assume that assumption (F.1) holds and $(Y^M,Z^M)$ is the solution of BSDE (\ref{EDSR approchee}).
\begin{itemize}
 \item Let assumption (B.2) holds. We have
$$e(M,n) \leqslant \frac{C}{e^{CM^2}}+\frac{Ce^{CM^{2rl}}}{n}.$$
In particular, for all $1<p < (rl)^{-1}$, if we take $M=(\log n)^{p/2}$ then $e(M,n)=o(h^{1-\varepsilon})$ for all $\varepsilon>0$.
 \item  Let assumption (B.1) holds. We have
$$e(M,n) \leqslant \frac{C}{e^{C_1M^2}}+\frac{Ce^{C_2M^{2}}}{n}.$$
In particular, if we take $M=\frac{1}{C_1+C_2}\sqrt{\log n}$ then $e(M,n)=o(h^{\frac{C_1}{C_1+C_2}})$.
\end{itemize}
\end{thm}

\paragraph{Proof of the theorem}
It is easy to see that 
$$e(M,n) \leqslant 2(e_1(M)+e_2(M,n))$$
with $e_1(M)$ defined by (\ref{erreur 1}) and
\begin{equation*}
 e_2(M,n) := \sup_{0\leqslant k\leqslant n} \E \left[ \abs{Y^{M,n}_{t_k}-Y^M_{t_k}}^2 \right] + \sum_{k=0}^{n-1} \E \left[ \int_{t_k}^{t_{k+1}} \abs{Z^{M,n}_{t_k}-Z^M_t}^2dt\right].
\end{equation*}
The error $e_1(M)$ is already studied in Proposition~\ref{proposition etude erreur 1}. Concerning the error $e_2(M,n)$ let us remark that $(Y^M,Z^M)$ is the solution of a BSDE with Lipschitz coefficients: indeed, $g_M$ and $f_M$ are Lipschitz with respect to $x$ and $y$, $f_M$ is locally Lipschitz with respect to $z$ and Proposition~\ref{prop Z borne} gives us that $Z^M$ is bounded. Thus, we are allowed to apply Proposition~\ref{prop resultat classique discretisation temporelle EDSR}:
$$e_2(M,n) \leqslant \frac{Ce^{CK_{f_M,z}^2}\left(1+K_{g_M}^2+K_{f_M,x}^2+\E\left[\int_0^T \abs{Z^M_t}^2dt \right]+\E \left[ \sup_{0 \leqslant t \leqslant T} \abs{Y^M_t}^2\right]\right)}{n}$$
with $K_{g_M}$ the Lipschitz constant of $g_M$, and $K_{f_M,x}$, $K_{f_M,z}$ the Lipschitz constants of $f_M$ with respect to $x$ and $z$. Estimations on $Z^M$ given by Proposition~\ref{proposition r<1/l} and Theorem~\ref{theoreme r=1surl} show us that $\E\left[\int_0^T \abs{Z^M_t}^2dt \right]$ is bounded by a bound that does not depend on $M$. Thanks to Itô's formula applied to $e^{K_{f_M,y}t}\abs{Y_t^M}^2$ and estimations on $Z^M$ given by Proposition~\ref{proposition r<1/l} and Theorem~\ref{theoreme r=1surl} it is also possible to show that $\E \left[ \sup_{0 \leqslant t \leqslant T} \abs{Y^M_t}^2\right]$ is also bounded by a bound that does not depend on $M$. Thus, we have
$$e_2(M,n) \leqslant \frac{Ce^{CK_{f_M,z}^2}\left(1+K_{g_M}^2+K_{f_M,x}^2\right)}{n}.$$
Under assumptions (F.1) and (B.2) we have, thanks to Proposition~\ref{prop Z borne}, $K_{g_M} \leqslant C(1+M^r)$, $K_{f_M,x} \leqslant C(1+M^r)$ and $K_{f_M,z} \leqslant C(1+(C(1+M^r))^l)\leqslant C(1+M^{rl})$. Finally, we obtain
$$e_2(M,n) \leqslant \frac{Ce^{CM^{2rl}}}{n}.$$
Under assumptions (F.1) and (B.1) we have, thanks to Proposition~\ref{prop Z borne}, $K_{g_M} \leqslant C(1+M^{1/l})$, $K_{f_M,x} \leqslant C(1+M^{1/l})$ and $K_{f_M,z} \leqslant C(1+(C(1+M^{1/l}))^l)\leqslant C(1+M)$. Finally, we obtain
$$e_2(M,n) \leqslant \frac{Ce^{CM^{2}}}{n}.$$
\cqfd

\begin{rem}
 Since $\sigma$ is a deterministic function, Euler and Milshtein schemes are equal, so the discretization error on $X$ is better. In this situation, authors of \cite{Gobet-Labart-07} show that the discretization error for the BSDE is on the same order than the discretization error for the SDE if we assume extra smoothness assumptions on $b$, $\sigma$, $g$ and $f$. More precisely, we could obtain the better estimate
 $$e_2(M,n) \leqslant \frac{Ce^{CM^{2}}}{n^2}.$$
\end{rem}

\begin{thm}
\label{theorem vitesse convergence sigma aleatoire}
We assume that assumptions (F.3) and (B.3) hold and $(\bar{Y}^M,\bar{Z}^M)$ is the solution of BSDE (\ref{EDSR approchee cas quadratique sigma aleatoire}).
\begin{itemize}
 \item If we assume that $2\kappa < 1-r$ then there exists $\eta$ such that
$$\bar{e}(M,n) \leqslant \frac{C}{e^{CM^{2+\eta}}}+\frac{Ce^{CM^2}}{n}.$$
In particular, for all $(2+\eta)^{-1}< p < 1/2$, if we take $M=(\log n)^{p}$ then $e(M,n)=o(h^{1-\varepsilon})$ for all $\varepsilon>0$.
 \item If we assume that $2\kappa = 1-r$ then we have
$$\bar{e}(M,n) \leqslant \frac{C}{e^{C_1M^2}}+\frac{Ce^{C_2M^{2}}}{n}.$$
In particular, if we take $M=\frac{1}{C_1+C_2}\sqrt{\log n}$ then $e(M,n)=o(h^{\frac{C_1}{C_1+C_2}})$.
\end{itemize}
\end{thm}

\paragraph{Proof of the theorem}
It is easy to see that 
$$\bar{e}(M,n) \leqslant 2(\bar{e}_1(M)+\bar{e}_2(M,n))$$
with $\bar{e}_1(M)$ defined by (\ref{erreur 1 bar}) and
\begin{equation*}
 \bar{e}_2(M,n) := \sup_{0\leqslant k\leqslant n} \E \left[ \abs{\bar{Y}^{M,n}_{t_k}-\bar{Y}^M_{t_k}}^2 \right] + \sum_{k=0}^{n-1} \E \left[ \int_{t_k}^{t_{k+1}} \abs{\bar{Z}^{M,n}_{t_k}-\bar{Z}^M_t}^2dt\right].
\end{equation*}
The error $\bar{e}_1(M)$ is already studied in Proposition~\ref{proposition etude erreur 1 bar}. Concerning the error $\bar{e}_2(M,n)$ let us remark that $(\bar{Y}^M,\bar{Z}^M)$ is the solution of a BSDE with Lipschitz coefficients: indeed, $g(\rho_{M^{(1-\kappa)^{-1}}}(.))$ and $f(.,\rho_{M^{(1-\kappa)^{-1}}}(.),.,\rho_{M}(.))$ are Lipschitz with respect to $x$, $y$ and $z$. Thus, we are allowed to apply Proposition~\ref{prop resultat classique discretisation temporelle EDSR}:
$$\bar{e}_2(M,n) \leqslant \frac{Ce^{CK_{f,z}^2}\left(1+K_g^2+K_{f,x}^2+\E\left[\int_0^T \abs{\bar{Z}^M_t}^2dt \right]+\E \left[ \sup_{0 \leqslant t \leqslant T} \abs{\bar{Y}^M_t}^2\right]\right)}{n}$$
with $K_g$ the Lipschitz constant of $g(\rho_{M^{(1-\kappa)^{-1}}}(.))$, and $K_{f,x}$, $K_{f,z}$ the Lipschitz constants of $f(.,\rho_{M^{(1-\kappa)^{-1}}}(.),.,\rho_{M}(.))$ with respect to $x$ and $z$. Classical estimates on solutions of quadratic BSDEs show us that $\E \left[ \sup_{0 \leqslant t \leqslant T} \abs{\bar{Y}^M_t}^2\right]$ and $\E\left[\int_0^T \abs{\bar{Z}^M_t}^2dt \right]$ are bounded by a bound that does not depend on $M$. Thus, we have
$$\bar{e}_2(M,n) \leqslant \frac{Ce^{CK_{f,z}^2}\left(1+K_g^2+K_{f,x}^2\right)}{n}.$$
Under assumption (B.3) we have, $K_g \leqslant C(1+M^{r(1-\kappa)})$, $K_{f,x} \leqslant C(1+M^{r(1-\kappa)})$ and $K_{f,z} \leqslant C(1+M)$. Finally, we obtain
$$\bar{e}_2(M,n) \leqslant \frac{Ce^{CM^{2}}}{n}.$$
\cqfd

\begin{rem}
 When $2\kappa >1-r$, the error estimate for $\bar{e}_1(M)$ given in remark \ref{remarque 2k >1-r} is not sufficient to obtain a ``good'' speed of convergence: the estimate on $\bar{e}(M,n)$ becomes, for $M$ well chosen,
$$\bar{e}(M,n) \leqslant \frac{C}{(\log n)^k},$$
for all $k \in \N^*$. This phenomenon is already explained in introductions of articles \cite{Richou-11,Imkeller-dosReis-09}.
\end{rem}

\appendix
\section{Appendix}

\subsection{Proof of Proposition~\ref{prop Z borne}}
To show the result we will use a classical truncation argument (see e.g. the proof of Theorem 4.1 in \cite{Delbaen-Hu-Bao-09}). Our truncation function $\rho_N$ is the projection on the centered euclidean ball of radius $N$ in $\R^{1 \times d}$. We denote $(Y^N,Z^N)$ the solution of the BSDE
$$Y_t^N=g(X_T)+\int_t^T f(s,X_s,Y_s^N,\rho_N(Z_s^N))ds-\int_t^T Z_s^NdW_s.$$
Now, this BSDE is also Lipschitz with respect to $z$. By the same calculus than in the proof of Theorem~3.1 in \cite{Richou-11} we can show that $Z^N$ is bounded by
$$\abs{Z^N} \leqslant e^{(2K_b+K_{f,y})T}\abs{\sigma}_{\infty}(K_g+TK_{f,x}).$$
This bound does not depend on $N$ so $\rho_N(Z^N)=Z^N$ for $N$ big enough. Then a uniqueness result for BSDEs with Lipschitz coefficients gives us that $(Y,Z)=(Y^N,Z^N)$ and the result is proved. \cqfd


\subsection{Proof of Lemma~\ref{lemme moment expo sous Q}}
Thanks to the estimate on $Z$ of Theorem~\ref{theoreme r=1surl} we easily show
\begin{eqnarray*}
 X_s &=&X_t+ \int_t^s b(u,X_u) du+\int_t^s \sigma(u)[d\tilde{W}_u+V_u^Mdu]\\
\abs{X_s}& \leqslant & \abs{X_t} + C + C\int_t^s\abs{X_u}du+\abs{\int_t^s\sigma(u)d\tilde{W}_u}\\
\sup_{t \leqslant r \leqslant s} \abs{X_r} &\leqslant& \abs{X_t}+C+C\int_t^s \sup_{t \leqslant r \leqslant u} \abs{X_r} du +\sup_{t \leqslant r \leqslant T} \abs{\int_t^r \sigma(u)d\tilde{W}_u},
\end{eqnarray*}
and we deduce from Gronwall's lemma the inequality
\begin{eqnarray*}
 \sup_{t \leqslant r \leqslant s} \abs{X_r} &\leqslant& C\left(1+\abs{X_t}+\sup_{t \leqslant r \leqslant T} \abs{\int_t^r \sigma(u)d\tilde{W}_u}\right).
\end{eqnarray*}
The first part of the lemma is easily proved thanks to the previous inequality. Moreover, we also have 
\begin{eqnarray*}
 \E_t^{\Q^M}\left[ e^{\mu \sup_{t \leqslant r \leqslant s} \abs{X_r}^2}\right] & \leqslant& Ce^{C\mu\abs{X_t}}\E_t^{\Q^M}\left[\sup_{t \leqslant s \leqslant T} \exp \left(C\mu\abs{\int_t^s \sigma(u)d\tilde{W}_u}^2\right)\right].
\end{eqnarray*}
It follows from the Dambis-Dubins-Schwarz representation theorem and Doob's maximal inequality that
$$\E_t^{\Q^M}\left[\sup_{t \leqslant s \leqslant T} \exp \left(C\mu\abs{\int_t^s \sigma(u)d\tilde{W}_u}^2\right)\right] \leqslant \E \left[ \sup_{0 \leqslant s \leqslant \abs{\sigma}_{\infty}^2(T-t)} e^{C\mu\abs{W_s}^2}\right] \leqslant 4\E \left[e^{C\mu\abs{\sigma}_{\infty}^2 T \abs{W_1}^2 }\right],$$
which is a finite constant if $C\mu\abs{\sigma}_{\infty}^2 T<1/2$. 
\cqfd

\subsection{Proof of Lemma~\ref{moment expo X sigma non borne}}
Let us consider the process
$$Y_t := \left(1+\abs{X_t}^2 \right)^{\frac{1-\kappa}{2}} = F(X_t).$$
Itô's formula gives us
\begin{eqnarray*}
 Y_t &=& Y_0 + (1-\kappa) \int_0^t \frac{^tX_sb(s,X_s)}{(1+\abs{X_s}^2)^{\frac{1+\kappa}{2}}}ds + \frac{1}{2}\int_0^t \mathrm{trace}\left(\nabla^2 F(X_s)\sigma(s,X_s)\tr{\sigma}(s,X_s)\right)ds\\
&& +(1-\kappa) \int_0^t \frac{\tr{X_s}\sigma(s,X_s)}{(1+\abs{X_s}^2)^{\frac{1+\kappa}{2}}}dW_s\\
 &=& Y_0+\int_0^t \tilde{b}(s,X_s)ds+\int_0^t \tilde{\sigma}(s,X_s) dW_s,\\
\end{eqnarray*}
with $\abs{\tilde{\sigma}}\leqslant C$ and $\abs{\tilde{b}(t,x)}\leqslant C(1+\abs{x}^{1-\kappa})$. Then, we have
\begin{eqnarray*}
 \sup_{0 \leqslant t \leqslant u} \abs{Y_t} &\leqslant& \abs{Y_0}+\int_0^u \sup_{0 \leqslant t \leqslant s} \abs{\tilde{b}(t,X_t)}ds + \sup_{0 \leqslant t \leqslant u} \abs{\int_0^t \tilde{\sigma}(s,X_s)dW_s}\\
&\leqslant& \abs{Y_0}+CT+ C\int_0^u \sup_{0 \leqslant t \leqslant s} \abs{X_t}^{1-\kappa}ds + \sup_{0 \leqslant t \leqslant T} \abs{\int_0^t \tilde{\sigma}(s,X_s)dW_s}\\
&\leqslant& \abs{Y_0}+CT+ C\int_0^u \sup_{0 \leqslant t \leqslant s} \abs{Y_t}ds + \sup_{0 \leqslant t \leqslant T} \abs{\int_0^t \tilde{\sigma}(s,X_s)dW_s},
\end{eqnarray*}
and we deduce from Gronwall's lemma the inequality
$$\sup_{0 \leqslant t \leqslant T} \abs{Y_t} \leqslant C\left( 1+ \sup_{0 \leqslant t \leqslant T} \abs{\int_0^t \tilde{\sigma}(s,X_s)dW_s} \right).$$
Since $\tilde{\sigma}$ is bounded, we are able to fit the end of Lemma~\ref{lemme moment expo sous Q} to show that there exists $\lambda>0$ such that
$$E\left[\exp \left( \lambda\sup_{0 \leqslant t \leqslant T} \abs{Y_t}^2 \right) \right] <+\infty.$$
Since $\abs{X_t}^{(1-\kappa)}\leqslant \abs{Y_t}$, the proof is complete. \cqfd

\subsection{Proof of Proposition~\ref{prop resultat classique discretisation temporelle EDSR}}
It is already proved in \cite{Bouchard-Touzi-04,Zhang-04} for the implicit scheme or in \cite{Lemor-05} for the explicit scheme that $e(n) =O(h)$. We just have to rewrite the proof to show where constants $K_g$, $K_{f,x}$ and $K_{f,z}$ appear precisely. For the readability of this article we will only give few key steps. Firstly, for the error in $Y$ we find, for $h$ small enough,
\begin{eqnarray}
 \nonumber \sup_{0\leqslant k \leqslant n} \E \abs{Y_{t_k}-Y_{t_k}^n}^2 &\leqslant& Ce^{CK_{f,z}^2} \Bigg[ (1+K_{f,z}^2)h\E\left[\sup_{0 \leqslant t \leqslant T} \abs{Y_t}^2 \right] + (1+K_{f,z}^4)h\E\left[\int_0^T \abs{Z_t}^2 dt \right]\\
\label{estimee 1}
&& +CK_{f,z}^2\sum_{k=0}^{n-1} \E \left[ \int_{t_k}^{t_{k+1}} \abs{Z_t -\bar{Z}_{t_k}}^2dt\right] \Bigg],
\end{eqnarray}
with $\bar{Z}_{t_k} :=\frac{1}{h} \E_{t_k} \int_{t_k}^{t_{k+1}} Z_t dt$. For the error in $Z$ we find, for $h$ small enough,
\begin{eqnarray}
\nonumber
\sum_{k=0}^{n-1} \E \int_{t_k}^{t_{k+1}} \abs{Z_t-Z_{t_k}^n}^2dt &\leqslant & Ch\left(1+K_{f,x}^2+K_g^2+\E\left[\sup_{0 \leqslant t \leqslant T} \abs{Y_t}^2 \right]+K_{f,z}^2\E\left[\int_0^T \abs{Z_t}^2 dt \right]\right)\\
\label{estimee 2}
&&+C\sum_{k=0}^{n-1} \E \left[ \int_{t_k}^{t_{k+1}} \abs{Z_t -\bar{Z}_{t_k}}^2dt\right]+CK_{f,z}^2 \sup_{0\leqslant k \leqslant n} \E \abs{Y_{t_k}-Y_{t_k}^n}^2.
\end{eqnarray}
The study of the error $\sum_{k=0}^{n-1} \E \left[ \int_{t_k}^{t_{k+1}} \abs{Z_t -\bar{Z}_{t_k}}^2dt\right]$ was done by Zhang in \cite{Zhang-04}. Theorem 3.5 in \cite{Richou-11} improve a little bit the estimate by studying how $K_g$ appears in the constant. Let us rewrite the proof of this theorem. We suppose in a first time that $b$, $\sigma$, $g$ and $f$ are differentiable with respect to $x$, $y$ and $z$. Then $Y$ and $Z$ are differentiable with respect to $x$ and we obtain that
$$ \sum_{k=0}^{n-1} \E \left[ \int_{t_k}^{t_{k+1}} \abs{Z_t -\bar{Z}_{t_k}}^2dt\right] \leqslant Ch\left(K_g^2+K_{f,x}^2+(1+K_{f,z}^2)\E \left[\int_0^T \abs{\nabla Z_t}^2dt\right]\right).$$
Thanks to classical estimates onto the solution of the BSDE solved by $(\nabla Y,\nabla Z)$ we have
$$\E \left[\int_0^T \abs{\nabla Z_s}^2ds\right] \leqslant C(1+K_g^2+K_{f,x}^2)(1+K_{f,z}^2).$$
Thus, we obtain 
\begin{equation}
\label{estimee 3}
\sum_{k=0}^{n-1} \E \left[ \int_{t_k}^{t_{k+1}} \abs{Z_t -\bar{Z}_{t_k}}^2dt\right] \leqslant Ch(1+K_g^2+K_{f,x}^2)(1+K_{f,z}^4).
\end{equation}
By standard approximation and stability results for Lipschitz BSDEs this estimate stays true when $b$, $\sigma$, $g$ and $f$ are not differentiable. Finally, by putting together (\ref{estimee 1}), (\ref{estimee 2}) and (\ref{estimee 3}), we have
$$e(n) \leqslant Che^{CK_{f,z}^2}\left[(1+K_{f,z}^8)(1+K_g^2+K_{f,x}^2)+(1+K_{f,z}^6)\E\left[\int_0^T \abs{Z_t}^2dt \right]+(1+K_{f,z}^4)\E \left[ \sup_{0 \leqslant t \leqslant T} \abs{Y_t}^2\right]\right],$$
and the final result can be easily deduced. \cqfd

\def\cprime{$'$}

\end{document}